\documentclass[reqno]{amsart}

\usepackage{amsmath}
\usepackage{amssymb}

\newcommand{\fs}{\triangle}
\newcommand{\bs}{\nabla}
\newcommand{\R}{\mathbb{R}}
\newcommand{\C}{\mathbb{C}}

\newcommand{\N}{\mathbb{N}}
\newcommand{\fa}{\mathfrak{a}}
\newcommand{\fc}{\mathfrak{c}}
\newcommand{\pd}{\partial}

\def\K2{\operatorname{K}}
\def\M2{\operatorname{M}}

\newcommand{\Id}{\mathcal{I}}
\newcommand{\I}{\mathfrak{I}}

\newtheorem{Theorem}{Theorem}[section]
\newtheorem{Lemma}[Theorem]{Lemma}
\newtheorem{Proposition}[Theorem]{Proposition}
\newtheorem{Corollary}[Theorem]{Corollary}

\theoremstyle{definition}
\newtheorem{Definition}[Theorem]{Definition}
\newtheorem{Remark}[Theorem]{Remark}
\newtheorem{Example}[Theorem]{Example}

\newtheorem*{Ch}{Charlier polynomials}
\newtheorem*{Kr}{Krawtchouk polynomials}
\newtheorem*{Me}{Meixner polynomials}
\newtheorem*{Ha}{Hahn polynomials}
\newtheorem*{HaT}{Hahn-type polynomials on $\N_0$}

\newcommand{\thref}[1]{Theorem \ref{#1}}

\newcommand{\prref}[1]{Proposition \ref{#1}}
\newcommand{\reref}[1]{Remark \ref{#1}}
\newcommand{\deref}[1]{Definition \ref{#1}}
\newcommand{\exref}[1]{Example \ref{#1}}
\newcommand{\coref}[1]{Corollary \ref{#1}}

\newcommand{\wt}{\widetilde}

\numberwithin{equation}{section}

\begin{document}

\title[Discrete orthogonal polynomials and difference equations]
{Discrete orthogonal polynomials and difference equations of several
variables}

\date{October 21, 2006}

\author{Plamen~Iliev}
\address{P.~Iliev, School of Mathematics, Georgia Institute of Technology, 
Atlanta, GA 30332--0160, USA}
\email{iliev@math.gatech.edu}

\author{Yuan~Xu}
\address{Y.~Xu, Department of Mathematics, University of Oregon, Eugene, 
OR 97403--1222, USA}
\email{yuan@math.uoregon.edu}


\begin{abstract}
The goal of this work is to characterize all second order difference 
operators of several variables that have discrete orthogonal polynomials as 
eigenfunctions. Under some mild assumptions, we give a complete solution of 
the problem.
\end{abstract}

\maketitle

\tableofcontents

\section{Introduction}\label{se1}

The goal of this study is to characterize all second order difference 
equations of several variables that have discrete orthogonal polynomials 
as eigenfunctions. More precisely, we consider difference operators 
of the form 
$$
D = \sum_{1\leq i, j\leq d}A_{i,j}\fs_i\bs_j  +\sum_{i=1}^dB_i\fs_i  +C \Id,
$$
where $\fs_i$ and $\bs_i$ are the forward and backward operator in the 
direction of the $i$-th coordinate of $\R^d$, respectively, $\Id$ is the 
identity operator,  $A_{i,j}$, $B_i$ and $C$ are functions of $x \in \R^d$.  
The discrete orthogonal polynomials in our consideration are polynomials that 
are orthogonal with respect to an inner product of the form
$$
   \langle f, g \rangle = \sum_{x\in V} f(x) g(x) W(x),
$$
where $V$ is a lattice set in $\R^d$ and $W$ is some positive weight 
function on $V$. There is a close correlation between $D$, $W$ and $V$.
Some restrictions need to be imposed on $V$ due to the complexity of the 
geometry in higher dimensions. Under some mild and, we believe, 
reasonable assumptions on $V$, we give a complete solution of the 
problem. 

For $d =1$, the one dimensional case, the classification problem was 
studied by several authors early in the last century, we refer to \cite{Al,NSU}
for references.  It was found that the classical discrete orthogonal 
polynomials, namely,  Hahn polynomials, Meixner polynomials, 
Krawtchouk polynomials, and Charlier polynomials are eigenfunctions 
of second order difference operators on the real line, and these are 
believed to be the only ones that have such property.  To our great 
surprise, however, another family of solutions turns up when we analyze 
the problem carefully.  The difference equation satisfied by this family 
of solutions is similar to that of Hahn polynomials but the parameters 
need to be chosen in a different way. In other words, the difference 
equation has two separate families of solutions when the parameters are
chosen differently.

The same phenomenon also appears in the case of two variables. The second 
order difference equations that have orthogonal polynomials as 
eigenfunctions were identified in \cite{Xu2}, but not all families of 
orthogonal polynomials were found. In fact, viewing it as an analogue of 
second order differential operators that have orthogonal polynomials as 
eigenfunctions, only one family of solutions that had been known in 
the literature was identified for each difference equation. The solutions
of difference equations, however, turn out to be far richer than that
of differential equations. For example, in the case of quadratic 
eigenvalues, only Hahn polynomials on the set $V= \{(x,y): x \ge 0, y \ge 0, 
x+y \le N\}$ are identified in \cite{Xu2}; there are in fact several other 
families, including orthogonal polynomials on $\N_0^2$, on $V = \{(x,y): 
0 \le x \le N_1, 0 \le y \le N_2\}$, and several others. This will be shown 
in the present paper.

We start with identifying difference operators $D$ that are self-adjoint
with respect to the inner product $\langle\cdot, \cdot \rangle$. This 
leads to compatibility conditions between the coefficients of the operator
$D$. The requirement that the difference equations have polynomial
solutions and some mild restrictions on $V$ reduce the coefficients of
$D$ to special simple forms, which allows us to use the compatibility
conditions to determine the coefficients of the difference operator.
A careful analysis of the result determines all possible solutions. 
It is well known that orthogonal bases in several variables
are not unique \cite{DX}.  For each solution, we give a family of 
mutually orthogonal polynomials explicitly and compute their norm;
in other words, we give a family of orthonormal basis explicitly. 

It should be mentioned that the analysis in several variables is by 
no means a straightforward extension of analysis in one or two variables. 
The complexity of the problem increases substantially as the 
dimension grows. For example, the number of solutions grows 
exponentially with the dimension, and the geometry of the admissible 
lattice sets becomes more involved as the dimension increases. 
Although our approach resembles the one used in \cite{Xu2}, several
new ideas are introduced to resolve various outstanding problems. 
The study in this paper is also much more systematic and thorough, 
as demonstrated by the new results obtained even in the case of $d=2$ 
and by the new family of discrete orthogonal polynomials of one variable.

The paper is organized as follows. The self-adjoint operators are 
treated in Section 2, where the compatibility conditions are
derived. In Section 3, we study difference equations that have 
polynomial solutions, called admissible equations, and show that 
the eigenvalues of such equations are either a quadratic polynomial or 
a linear polynomial of the index, and the coefficients of such equations
need to have certain simple forms under some mild restrictions on $V$.
The compatibility conditions can then be used to determine the 
coefficients, thus $D$,  explicitly. The case of one variable is discussed 
in Section 4, in which the new family of discrete orthogonal polynomials 
is treated in details. For several variables, the case of  quadratic
eigenvalues is studied in Section 5 and the case of linear eigenvalues 
is developed in Section 6.

\section{Self-adjoint difference operators}\label{se2}
In this section we discuss the 
basic properties of self-adjoint second order difference operators with 
respect to an inner product defined on a discrete set. We start with some 
notations. 

Let $V$ be a subset of $\R^d$ and 
let $W:V\rightarrow \R$ be a positive function. 
Consider the inner product defined by the weight $W$ as follows
\begin{equation}\label{2.1}
\langle f,g\rangle=\sum_{x\in V} f(x)g(x)W(x),
\end{equation}
for functions $f,g:\R^d\rightarrow\R$. We will be mainly interested in 
the case when $f(x)$ and $g(x)$ are polynomials of $x\in\R^d$, but the 
general discussion in this section does not depend essentially on the class 
of functions. Denote by $\{e_1,e_2,\dots,e_d\}$ the standard basis for $\R^d$. 
We denote by $E_i$, $\fs_i$ and $\bs_i$, respectively,  the customary 
shift, forward and backward difference operators acting on a function $f(x)$ 
as follows
\begin{align*}
&E_if(x)=f(x+e_i)\\
&\fs_i f(x)=f(x+e_i)-f(x)=(E_i-\Id)f(x)\\
&\bs_i f(x)=f(x)-f(x-e_i)=(\Id-E_i^{-1})f(x),
\end{align*}
where $\Id$ is the identity operator.

We will consider second order difference operators of the form
\begin{equation}\label{2.2}
D=\sum_{1\leq i, j\leq d}A_{i,j}\fs_i\bs_j+\sum_{i=1}^dB_i\fs_i+C\Id,
\end{equation}
where $A_{i,j}$, $B_i$ and $C$ are some functions of $x$. Immediately 
from the definition one can see that $D$ can also be rewritten as
\begin{equation}\label{2.3}
D=\sum_{1\leq i\neq j\leq d}\alpha_{i,j}E_iE_j^{-1}+\sum_{i=1}^d\beta_iE_i+
\sum_{i=1}^d\gamma_iE_i^{-1}+\delta\Id,
\end{equation}
where the new coefficients $\alpha_{i,j}$, $\beta_i$, $\gamma_i$ and $\delta$ 
are related to the old ones via the formulas
\begin{align}
&\alpha_{i,j}=-A_{i,j} \text{ for }1\leq i\neq j\leq d\label{2.4}\\
&\beta_{i}=\sum_{k=1}^{d}A_{i,k}+B_i \text{ for }1\leq i\leq d\label{2.5}\\
&\gamma_{i}=\sum_{k=1}^dA_{k,i}\label{2.6}\\
&\delta=C-\sum_{1\leq i\neq j\leq d}A_{i,j}-\sum_{i=1}^d(2A_{i,i}+B_i).
\label{2.7}
\end{align}
The representation \eqref{2.2} will be more convenient when we deal with 
polynomials, because the operators $\fs_i$ and $\bs_i$ decrease the total 
degree of a polynomial by one. However, necessary and sufficient conditions 
for an operator $D$ to be self-adjoint with respect to the inner product 
\eqref{2.1} are much simpler and natural if we write it as in \eqref{2.3}. 
Formulas \eqref{2.4}-\eqref{2.7} allow us to go easily  from one 
representation to another.

Let us define ``directional'' boundaries of $V$ as follows:
\begin{align*}
&\pd_j^{\pm} V=\{x\in V: x\pm e_j\notin V\}&&\text{ for }j=1,2,\dots,d\\
&\pd_{i,j} V=\{x\in V: x+e_i- e_j\notin V\}&&\text{ for }1\leq i\neq j\leq d.
\end{align*}

\begin{Example}\label{ex2.1}
Let $\N_0$ denote the set of all nonnegative integers. 
For the set $V=\N_0^d$ we have $\pd_j^{-}V=\pd_{i,j}V=V\cap\{x:x_j=0\}$ and 
$\pd_{j}^{+}V=\emptyset$ for $j=1,2,\dots,d$ and $i\neq j$.

Other sets that play crucial role later are listed below
\begin{subequations}\label{2.8}
\begin{align}
&V^d_l=\{x\in\N_0^d:x_i\leq l_i\text{ for }i=1,2,\dots,d\}  \label{2.8a},\\
\intertext{where $l_i$ are positive integers;}
&V^d_N=\{x\in \N_0^d:x_1+x_2+\cdots+x_d\leq N\},  \label{2.8b}
\intertext{where $N$ is a positive integer;}
&V^d_{N,S}=V^d_N \bigcap_{i\in S}\{x:x\leq l_i\text{ for }i\in S\},\label{2.8c}
\end{align}
\end{subequations}
where $S$ is a nonempty subset of $\{1,2,\dots, d\}$ and $l_i$ are integers 
such that $1\leq l_i\leq N$.

For example, for the parallelepiped $V^d_l$ we have
$\pd_j^{-}V^d_l=V^d_l\cap\{x:x_j=0\}$, 
$\pd_{j}^{+}V^d_l=V^d_l\cap\{x:x_j=l_j\}$ and 
$\pd_{i,j}V^d_l=\pd_j^{-}V^d_l\cup\pd_i^{+}V^d_l$.

For $V^d_N$ defined by \eqref{2.8b} we have
$\pd_j^{-}V^d_N=\pd_{i,j}V^d_N=V^d_l\cap\{x:x_j=0\}$, and 
$\pd_{j}^{+}V^d_N=V^d_l\cap\{x:x_1+x_2\cdots+x_d=N\}$.
\end{Example}

The next proposition characterizes the self-adjoint operators 
with respect to the inner product \eqref{2.1} in terms of their coefficients.

\begin{Proposition} \label{pr2.2}
The operator $D$ is self-adjoint with respect to the inner 
product \eqref{2.1}, if and only if
\begin{align}
& W(x) \gamma_i(x) = W(x-e_i) \beta_i(x-e_i), \qquad x, x-e_i \in V  
\label{2.9}\\ 
&  W(x-e_i) \alpha_{i,j}(x-e_i) = W(x-e_j) \alpha_{j,i}(x-e_j), 
\quad x-e_i, x-e_j \in V, \quad j \ne i, \label{2.10}
\end{align}
and
\begin{subequations}\label{2.11}
\begin{align}
\gamma_j(x) & = 0, \qquad x \in \pd_j^{-} V  \label{2.11a},\\ 
\beta_j(x) & = 0, \qquad x \in \pd_j^{+} V  \label{2.11b}, \\ 
\alpha_{i,j}(x) & = 0, \qquad  x \in \pd_{i,j} V, \quad i \ne j.\label{2.11c}
\end{align}
\end{subequations}
\end{Proposition}

\begin{proof}
Using 
$\{x+e_i-e_j:x\in V\setminus \pd_{i,j}V\}=V\setminus \pd_{j,i} V$
and changing the summation index we see that 
\begin{align*}
&\langle \alpha_{i,j}E_iE_j^{-1} u,v\rangle=
\sum_{x \in \pd_{i,j}V} 
    \alpha_{i,j}(x)u(x+e_i-e_j) v(x)W(x)\\
&\qquad +\sum_{x \in V\setminus\pd_{j,i}V} 
    \alpha_{i,j}(x+e_j-e_i) u(x)v(x+e_j-e_i)W(x+e_j-e_i).
\end{align*}
Writing similar equalities for $\langle\beta_iE_i u,v\rangle$ and 
$\langle\gamma_iE_i^{-1} u,v\rangle$, using 
$\{x\pm e_i:x\in V\setminus \pd_{i}^{\pm}V\}=V\setminus \pd_{i}^{\mp} V$,
and comparing $\langle D u, v\rangle$ with  $\langle u, D v\rangle$ 
gives the stated conditions.
\end{proof}

As an immediate corollary from the relations established above we 
can deduce compatibility conditions that need to be satisfied by the 
coefficients of the operator $D$.

\begin{Corollary}[Compatibility conditions] \label{co2.3}
If the operator $D$ is 
self-adjoint with respect to the inner product defined by \eqref{2.1} then for 
$1\leq i\neq j\leq d$ we have
\begin{subequations}\label{2.12}
\begin{equation}\label{2.12a}
\alpha_{i,j}(x-e_i)\beta_j(x-e_j)\gamma_i(x)=
\alpha_{j,i}(x-e_j)\beta_i(x-e_i)\gamma_j(x),
\end{equation}
for $x,x-e_i,x-e_j\in V$ and
\begin{align}
&\beta_i(x-e_i)\beta_j(x-e_i-e_j)\gamma_j(x)\gamma_i(x-e_j)=\nonumber\\
&\qquad\qquad\beta_j(x-e_j)\beta_i(x-e_i-e_j)\gamma_i(x)\gamma_j(x-e_i),
\label{2.12b}
\end{align}
for $x,x-e_i,x-e_j,x-e_i-e_j\in V$.
\end{subequations}
\end{Corollary}

\begin{proof} Writing 
\begin{equation*}
\frac{W(x-e_j)}{W(x-e_i)}=\frac{W(x-e_j)}{W(x)}\frac{W(x)}{W(x-e_i)}
\end{equation*}
and using \eqref{2.10} on the left side and \eqref{2.9} for the right side
we obtain \eqref{2.12a}. Similarly if we use \eqref{2.9} for all ratios in
the identity
\begin{equation*}
\frac{W(x)}{W(x-e_i)}\frac{W(x-e_i)}{W(x-e_i-e_j)}
=\frac{W(x)}{W(x-e_j)}\frac{W(x-e_j)}{W(x-e_i-e_j)},
\end{equation*}
we get \eqref{2.12b}, which completes the proof.
\end{proof}

\begin{Remark}\label{re2.4} 
The compatibility conditions stated above can be easily extended to
more general difference operators. For example, we can consider 
the operator $\wt D$ that adds two additional terms $E_iE_j$ and 
$E_i^{-1}E_j^{-1}$ in the operator $D$ in \eqref{2.3}. However, the 
complexity of the computations in the sections below will increase 
significantly even with these two terms added. Furthermore, since 
$\fs = \bs + \bs \fs$, the operator $\wt D$ includes the 4-th order term 
$\fs_1 \bs_1 \fs_2 \bs_2$ (say, $d =2$) which is the discrete analog of 
$\partial_1^2\partial_2^2$, where 
$\partial_i$ stands for the partial derivative with respect to $x_i$.
Classifying orthogonal polynomials that are eigenfunctions of $\wt D$ 
appears to be an interesting but much harder problem.
\end{Remark}

\section{Admissible equations and orthogonal polynomials}\label{se3}

In this section we study difference operators $D$ that have discrete 
orthogonal polynomials as eigenfunctions. For more details on the general 
theory of discrete orthogonal polynomials of several variables we refer the 
reader to \cite{Xu1}.

Throughout the paper we use the standard multi-index notation. A 
multi-index will be denoted by 
$\mu=(\mu_1,\mu_2,\dots,\mu_d)\in\N^d_0$. For each $\mu$ we denote by $x^\mu$
the monomial 
\begin{equation*}
x^{\mu}=x_1^{\mu_1}x_2^{\mu_2}\cdots x_d^{\mu_d}
\end{equation*}
of total degree $|\mu|=\mu_1+\mu_2+\cdots+\mu_d$. The degree of a polynomial 
is defined as the highest degree of its monomials. We denote by 
$\R[x]=\R[x_1,x_2,\dots,x_d]$ the space of all polynomials in the variables 
$x_1,x_2,\dots,x_d$ and by $\Pi^d_n$ the subspace of polynomials of degree at 
most $n$ in the variables $\{x_1,x_2,\dots,x_d\}$. The latter has dimension 
$\dim(\Pi^d_n)=\binom{n+d}{n}$.

Let $V$ be an at most countable set of isolated points in $\R^d$ and let 
$|V|$ denote the cardinality of $V$. If 
$\langle f,g \rangle=0$ for the inner product defined by \eqref{2.1}, 
we say that $f$ and $g$ are orthogonal with respect to $W$ on the set $V$. 
Orthogonal polynomials on $V$ depend on the structure of the ideal
\begin{equation*}
\I(V)=\{p\in\R[x_1,x_2,\dots,x_d]:p(x)=0\text{ for every }x\in V\}.
\end{equation*}
In fact, the orthogonal polynomials belong to the space
\begin{equation*}
\R[V]=\R[x_1,x_2,\dots,x_d]/\I(V).
\end{equation*}
There is a lattice set $\Lambda(V)$ such that every polynomial in 
$\R[V]$ can be written as
\begin{equation*}
P(x)=\sum_{\mu\in\Lambda(V)}c_{\mu}x^{\mu}.
\end{equation*}
The lattice set is not uniquely determined by $V$. However, in all cases 
that arise here, there is a natural way to identify $V$ and $\Lambda(V)$, 
see \exref{ex3.1}.
Denote $\Lambda_k(V)=\{\mu\in\Lambda(V):|\mu|=k\}$ and let $r_k$ be the 
number of the elements in the set $\Lambda_k(V)$. 

\begin{Example}\label{ex3.1} 
It is clear that for $V=\N_0^d$ we have $\I(\N_0^d)=(0)$ and therefore 
$\R[\N_0^d]=\R[x]$. For all other sets defined in \exref{ex2.1} one can 
easily find (inductively on $d$ and $|V|$) explicit sets of generators of the 
ideal $\I(V)$. In particular, the generators listed below allow us to 
identify $\Lambda(V)$ and $V$.

For $V^d_l$, defined by \eqref{2.8a}, the ideal $\I(V^d_l)$ is generated 
by the set
$$G^d_l=\{(-x_1)_{l_1+1},(-x_2)_{l_2+1},\dots,(-x_d)_{l_d+1}\}.$$ 

For $V^d_N$, given by \eqref{2.8b}, the ideal $\I(V^d_N)$ is generated by 
the set
$$G^d_N=\{(-x_1)_{\mu_1}(-x_2)_{\mu_2}\cdots(-x_d)_{\mu_d}:|\mu|=N+1\}.$$ 

Finally, for $V^d_{N,S}$ in \eqref{2.8c}, the corresponding ideal is 
generated by 
$$G^d_{N,S}=G^d_N\cup\{(-x_i)_{l_i+1}:i\in S\}.$$
\end{Example}

Next we want to consider difference operators $D$ on the space $\R[V]$. In 
order to have a proper action of $D$ on $\R[V]$ the ideal $\I(V)$ must be 
$D$-invariant, i.e. 
\begin{equation*}
D(\I(V))\subset \I(V).
\end{equation*}
One can easily show that if $D$ is self-adjoint with respect to the inner 
product \eqref{2.1}, then the above condition is satisfied. More precisely, 
we have the following proposition.

\begin{Proposition}\label{pr3.2}
If the operator $D$, written in the form \eqref{2.3}, satisfies conditions 
\eqref{2.11} then the ideal $\I(V)$ is $D$-invariant.
\end{Proposition}

\begin{proof}
Notice that if a polynomial $p(x)$ vanishes on $V$, then 
the polynomial $p(x-e_j)$ will vanish on $V\setminus \pd_j^{-}V$. Thus, if 
$p(x)\in \I(V)$ and \eqref{2.11a} holds then $\gamma_j(x)p(x-e_j)\in\I(V)$, 
i.e. the operator $\gamma_jE_j^{-1}$ preserves $\I(V)$. Similar arguments show 
that \eqref{2.11b} (resp. \eqref{2.11c}) implies that 
$\beta_jE_j$ (resp. $\alpha_{i,j}E_iE_j^{-1}$) preserves $\I(V)$, which 
completes the proof.
\end{proof}

\begin{Definition} \label{de3.3} 
Let $\I(V)$ be $D$-invariant. The equation 
\begin{equation}\label{3.1}
Du=\lambda u
\end{equation}
is called admissible on $V$ if for any $k\in\N_0$ there is a number 
$\lambda_k$ such that the equation $Du=\lambda_k u$ has $r_k$ linearly 
independent polynomial solutions in $\R[V]$ and it has no nontrivial 
solutions in the set of polynomials of degree less than $k$.
\end{Definition}

\begin{Remark}\label{re3.4}
Notice that if \eqref{3.1} is admissible on $V$ and if we have nontrivial 
solutions of the equations $Du=\lambda_k u$ and $Du=\lambda_l u$ for 
distinct integers $k,l\in\N_0$ then $\lambda_k\neq\lambda_l$. Thus if 
$D$ is also a self-adjoint operator with respect to the inner product 
\eqref{2.1}, the polynomial solutions for 
$Du=\lambda_k u$ and $Du=\lambda_l u$ will be mutually orthogonal. 
\end{Remark}

From now on, we will assume that $\lambda_0=0$ (which is 
equivalent to $C=0$ for operators $D$ of the form \eqref{2.2}). This 
is not a real restriction on $D$ since one can replace $D$ by $D-\lambda_0\Id$
if necessary.

The next proposition shows that we can pick linearly independent 
polynomial solutions in such a way that the highest terms contain 
single monomials.

\begin{Proposition} \label{pr3.5}
The equation \eqref{3.1} is admissible on $V$ if and only if for each 
$k\in\N_0$ there exists a number $\lambda_k$ such that the equation 
$Du=\lambda_k u$ has $r_k$ linearly independent polynomial solutions of 
the form
\begin{equation}
P_{\mu}(x)=x^{\mu} \mod \Pi^d_{k-1}.
\end{equation}
\end{Proposition}

The proof of \prref{pr3.5} follows easily form \deref{de3.3}, see \cite{Xu2} 
in the case $d=2$.

\begin{Remark}\label{re3.6}
If $D$ is an admissible operator on $V$ we can trivially extend it to 
an admissible operator on the set $V'=\{(x,0):x\in V\}\subset\R^{d+1}$.
When we classify the possible admissible operators, we will exclude 
such degenerate situations.
\end{Remark}

The next proposition essentially characterizes the admissible difference 
operators of the form \eqref{2.2}. In order to eliminate several singular  
cases we will assume that $\mu\in\Lambda(V)$ for $|\mu|=1,2,3$. Instead of 
putting $\mu\in\Lambda(V)$ for $|\mu|=3$ we can require that 
$3e_i\in\Lambda(V)$ for some $i\in\{1,2,\dots,d\}$ and for all 
$1\leq i\neq j\leq d$ at least one of $\{2e_i+e_j,e_i+2e_j\}$ belongs to 
$\Lambda(V)$, but this would make the statement of the theorem long 
and awkward.

\begin{Theorem}\label{th3.7}
Assume that $\mu\in\Lambda(V)$ for $|\mu|=1,2,3$. Let $D$
be a second order difference operator defined by \eqref{2.2} such that
$\I(V)$ is $D$-invariant. Then the following conditions are equivalent:
\begin{itemize}
\item[(i)] The equation $Du=\lambda u$ is admissible;
\item[(ii)] 
The coefficients $B_{i}(x)$ are polynomials of degree at most 1, and 
$A_{i,j}(x)$ are polynomials of degree at most 2 satisfying
\begin{subequations}\label{3.3}
\begin{align}
&B_{i}  = bx_i \mod \Pi^d_0 \text{ for }i=1,2,\dots,d\label{3.3a}\\
&A_{i,i}= ax_i^2\mod \Pi^d_1\text{ for }i=1,2,\dots,d\label{3.3b}\\
&A_{i,j}+A_{j,i}= 2ax_ix_j \mod \Pi^d_1 \text{ for }1\leq i\neq j\leq d
                                   \label{3.3c}\end{align}
\end{subequations}
for some constants $a$ and $b$. The eigenvalues $\lambda_k$ must be distinct 
and they are given by 
\begin{equation}\label{3.4}
\lambda_k=k(ka-a+b).
\end{equation}
\end{itemize}
\end{Theorem}

\begin{proof} Assume first that the equation $Du=\lambda u$ is admissible. 
Applying \prref{pr3.5} we see that there exist polynomials
\begin{equation*}
P_{e_i}(x)=x_i \mod \Pi^d_0 
\end{equation*}
satisfying 
\begin{equation}\label{3.5}
Du=\lambda_1u.
\end{equation}
Notice that $\fs_j P_{e_i}=\delta_{i,j}$ and $\fs_j\bs_kP_{e_i}=0$ for all 
$j,k=1,2,\dots,d$. Thus $DP_{e_i}=B_i(x)$ and the last equation shows that 
formula \eqref{3.3a} holds, where we put $b=\lambda_1$. 
Similarly, we can find polynomials of the form
\begin{equation}\label{3.6}
P_{2e_i}(x)=x_i^2 \mod \Pi^d_1 
\end{equation}
satisfying
\begin{equation}\label{3.7}
Du=\lambda_2u.
\end{equation}
This time we have 
$\fs_j P_{2e_i}=2x_i\delta_{i,j} \mod \Pi^d_0$ and 
$\fs_j\bs_kP_{2e_i}=2\delta_{i,j}\delta_{i,k}$.
Using now \eqref{3.3a} we see that
\begin{equation*}
DP_{2e_i}=2A_{i,i}+2bx_i^2\mod \Pi^d_1,
\end{equation*}
which, combined with \eqref{3.6} and \eqref{3.7}, shows that
\begin{equation*}
A_{i,i}=\frac{\lambda_2-2b}{2}x_i^2\mod \Pi^d_1,
\end{equation*}
thus proving \eqref{3.3b} if we denote $a=(\lambda_2-2b)/2$.

Next, for $i\neq j$ we take a polynomial of the form
\begin{equation*}
P_{e_i+e_j}(x)=x_ix_j \mod \Pi^d_1 
\end{equation*}
satisfying \eqref{3.7}. A simple computation as above gives that 
\begin{equation*}
DP_{e_i+e_j}=A_{i,j}+A_{j,i}+2bx_ix_j\mod \Pi^d_1,
\end{equation*}
which combined with \eqref{3.7} shows that \eqref{3.3c} holds. It remains to 
show that $A_{i,j}$ are polynomials of degree at most 2 and that \eqref{3.4} 
holds. First we use a polynomial of the form 
\begin{equation*}
P_{3e_i}=x_i^3 \mod \Pi^d_2
\end{equation*}
for some $i$ satisfying 
\begin{equation}\label{3.8}
Du=\lambda_3u,
\end{equation}
to show that \eqref{3.4} holds for $k=3$. Indeed, we already know that 
if $P$ is a polynomial of degree at most 2, $D(P)$ will also be a polynomial 
of degree at most 2 (because we can write it as a linear combination of $1$,
$P_{e_i}$, $P_{e_i+e_j}$, for $i,j=1,2,\dots,d$). Thus 
$D(P_{3e_i})=D(x_i^3)\mod \Pi^d_2$. On the other hand
\begin{equation*} 
D(x_i^3)=6x_iA_{i,i}+(3x_i^2+3x_i+1)B_{i}=(6a+3b)x_i^3\mod \Pi^d_2,
\end{equation*}
where in the last equality we used \eqref{3.3a} and \eqref{3.3b}. Comparing 
now the coefficients of $x_i^3$ on both sides in \eqref{3.8}, we see that 
$\lambda_3=6a+3b$, which is exactly \eqref{3.4} for $k=3$.

Next we consider the solution to equation \eqref{3.8} of the form
\begin{equation*}
P_{2e_i+e_j}(x)=x_i^2x_j \mod \Pi^d_2
\end{equation*}
for $i\neq j$. Again $D(P_{2e_i+e_j})=D(x_i^2x_j)\mod\Pi^d_2$, but this time 
\begin{equation*}
\begin{split} 
D(x_i^2x_j)
&=(2x_i+1)A_{i,j}+(2x_i-1)A_{j,i}+2x_jA_{i,i}+(2x_i+1)x_jB_i+x_i^2B_j\\
&=2x_i(A_{i,j}+A_{j,i})+A_{i,j}-A_{j,i}+2x_jA_{i,i}+2x_ix_jB_i+x_i^2B_j
                                                \mod \Pi^d_2\\
&=(6a+3b)x_i^2x_j+A_{i,j}-A_{j,i} \mod \Pi^d_2\\
&=\lambda_3 x_i^2x_j+A_{i,j}-A_{j,i}\mod \Pi^d_2,
\end{split}
\end{equation*}
upon using \eqref{3.3} and $\lambda_3=6a+3b$.
Equation \eqref{3.8} shows that $A_{i,j}-A_{j,i}=0\mod \Pi^d_2$, which 
combined with \eqref{3.3c} proves that $A_{i,j}$ and $A_{j,i}$ are polynomials 
of degree at most 2. Finally, let $\mu\in\Lambda(V)$ such that $|\mu|=k$ and 
let 
\begin{equation*}
P_{\mu}=x^{\mu} \mod \Pi^d_{k-1}
\end{equation*}
be a solution to
\begin{equation}\label{3.9}
Du=\lambda_ku.
\end{equation}
The admissibility implies that $D(u)=D(x^{\mu})\mod \Pi^d_{k-1}$. 
Notice that 
\begin{subequations}
\begin{align}
&\fs_i x^{\mu}
      =\mu_ix^{\mu-e_i}\mod \Pi^d_{k-2}\text{ for all }i\\
&\fs_i\bs_j x^{\mu}
      =\mu_i\mu_jx^{\mu-e_i-e_j}\mod \Pi^d_{k-3}\text{ for }i\neq j\\
&\fs_i\bs_i x^{\mu}
      =\mu_i(\mu_i-1)x^{\mu-2e_i}\mod \Pi^d_{k-3}\text{ for all }i
\end{align}
\end{subequations}
and therefore
\begin{equation}\label{3.11}
\begin{split} 
D(x^{\mu})
&=\sum_{1\leq i< j\leq d} \mu_i\mu_jx^{\mu-e_i-e_j}(A_{i,j}+A_{j,i})+
\sum_{i=1}^d \mu_i(\mu_i-1)x^{\mu-2e_i}A_{i,i}\\
&\qquad+\sum_{i=1}^d \mu_ix^{\mu-e_i}B_{i}\mod \Pi^d_{k-1}\\
&=x^{\mu}\left(2a\sum_{1\leq i< j\leq d} \mu_i\mu_j+
a\sum_{i=1}^d \mu_i(\mu_i-1)+b\sum_{i=1}^d \mu_i\right)\mod \Pi^d_{k-1}\\
&=|\mu|(a|\mu|-a+b)x^{\mu}\mod \Pi^d_{k-1},
\end{split}
\end{equation}
which combined with \eqref{3.9} gives \eqref{3.4}

Conversely, assume that (ii) holds. One can easily show by induction on 
$|\mu|$ that for $\mu\in\N^d$ there exist polynomials of the form
\begin{equation}\label{3.12}
P_{\mu}(x)=x^{\mu}+\sum_{|\nu|<|\mu|}c_{\mu,\nu}P_{\nu}(x)
\end{equation}
satisfying \eqref{3.9} with eigenvalue $\lambda_{|\mu|}$ given by 
\eqref{3.4}. Indeed, let us assume that this is true for $|\mu|\leq k-1$ and 
take $\mu$ such that $|\mu|=k$. Using computation \eqref{3.11} we see that 
\begin{equation*}
D(x^{\mu})=\lambda_{|\mu|}x^{\mu}
+\sum_{|\nu|<|\mu|}\gamma_{\mu,\nu}P_{\nu}(x),
\end{equation*}
for some constants $\gamma_{\mu,\nu}$. Then 
$$D(P_{\mu})-\lambda_{|\mu|} P_{\mu}=\sum_{|\nu|<|\mu|}
(c_{\mu,\nu}(\lambda_{|\nu|}-\lambda_{|\mu|})+\gamma_{\mu,\nu})P_{\nu},$$ 
i.e. if we pick 
$c_{\mu,\nu}=-\gamma_{\mu,\nu}/(\lambda_{|\nu|}-\lambda_{|\mu|})$ 
the polynomial $P_{\mu}$ defined by \eqref{3.12} will satisfy \eqref{3.9}.
Thus, in the quotient space $\R[V]$ we have $r_k$ polynomials of total 
degree $k$ satisfying \eqref{3.9}.
\end{proof}

In the next sections we will use conditions \eqref{3.3} and the 
compatibility conditions \eqref{2.11} and \eqref{2.12} to 
determine the possible self-adjoint operators on appropriate sets $V$. 

Since the difference equation $Du=\lambda u$ is invariant under translations 
$x\rightarrow x+h$, we can  consider the difference equations modulo
translations. 

In the rest of the paper we will impose certain boundary conditions on the 
difference operator $D$. More precisely, we will assume that for every 
$j$ the coefficients $A_{i,j}(x)$ vanish on the 
hyperplane $\{x:x_j=0\}$. Roughly speaking, this means that $V$ has 
sufficiently many boundary points belonging to $\{x:x_j=0\}$. We show 
below, that if $a\neq 0$ and if $V$ has ``enough'' points on its boundary 
then, up to a translation, this must be true.

Let us denote by $\pd_j V=\pd_j^{-}V \cap_{i\neq j}\pd_{i,j}V$ the ``lower'' 
$j$th boundary of $V$. We say that $V$ has a nontrivial  boundary 
if for every $j=1,2,\dots,d$, $\pd_j V$ contains more than $2^d$ points and 
that these points do not belong to a variety of dimension less than $d-1$.
For example, a line in $d=2$ is a nontrivial boundary, so is a plane in $d=3$. 

\begin{Proposition}\label{pr3.8}
Let $V$ be a discrete set with a nontrivial boundary and let $D$ be an 
admissible, self-adjoint operator. If $a\neq 0$, then after an appropriate 
translation, $\partial_j V$ becomes a subset of $\{x: x_j =0\}$ and the 
coefficients $A_{i,j}(x)$ vanish on $\{x: x_j =0\}$.
\end{Proposition}

\begin{proof}
Using the fact that $D$ is self-adjoint, \prref{pr2.2}, \eqref{2.4} 
and \eqref{2.6} we see that 
$$A_{i,j}(x)   = 0, \qquad x \in \pd_j V \text{ for } 1 \le j \le d.$$
Applying Bezout's theorem for polynomials in $d$ variables 
\cite[Chapt. IV, Section 2]{Sh} for the polynomials 
$A_{1,j}(x),\dots, A_{d,j}(x)$ which satisfy \eqref{3.3b}-\eqref{3.3c} and 
vanish on $\pd_j V$, one concludes that they must contain a common linear 
factor. Equation \eqref{3.3b} shows that this factor must be of the form 
$x_j-h_j$. Thus, applying the translation $x\rightarrow x+h$ with 
$h=(h_1,h_2,\dots,h_d)$, $A_{i,j}$ takes the form
\begin{equation}\label{3.13}
A_{i,j}(x) = x_j\times(\text{linear polynomial in }x_1,x_2,\dots,x_d),
\end{equation}
and completes the proof. 
 \end{proof}

The compatibility conditions also imply that $A_{i,j}-A_{j,i}\in\Pi^d_1$ for 
$i\neq j$, which combined with \eqref{3.3c} shows that 
$A_{i,j}=ax_ix_j\mod\Pi^d_1$ (the argument works simultaneously for 
$a\neq 0$ and $a=0$). The form of the coefficients determined below will be 
the staring point in the next sections.

\begin{Proposition}\label{pr3.9}
Let $A_{i,j}$ and $B_{i}$ satisfy the conditions in \thref{th3.7}(ii) and 
assume that $A_{i,j}$ vanishes on $\{x:x_j=0\}$. Then, the polynomial 
identities \eqref{2.12a} imply that
\begin{align}
&A_{i,j}=x_j(ax_i+l_{i,j})\label{3.14}\\
&B_i=bx_i+s_i,\label{3.15}
\end{align}
or equivalently:
\begin{align}
&\alpha_{i,j}=-x_j(ax_i+l_{i,j})\label{3.16}\\
&\beta_i=ax_i\sum_{k=1}^dx_k+\sum_{k=1}^dl_{i,k}x_k+bx_i+s_i\label{3.17}\\
&\gamma_i=ax_i\sum_{k=1}^dx_k+x_i\sum_{k=1}^dl_{k,i},\label{3.18}
\end{align}
for some constants $a$, $b$, $l_{i,j}$ and $s_i$.
\end{Proposition}

\begin{proof}
Equation \eqref{3.15} follows immediately from \eqref{3.3a}. 
From equations \eqref{3.3b}, \eqref{3.3c} and the fact that $A_{i,j}$ vanishes 
on $\{x:x_j=0\}$ we deduce that
\begin{equation}
A_{i,j}=x_j(q_{i,j}x_i+l_{i,j}),
\end{equation}
where $q_{i,i}=a$ and $q_{i,j}+q_{j,i}=2a$. 
It remains to show that $q_{i,j}=a$ for all $i\neq j$. 
From \eqref{2.4}-\eqref{2.6} we see that 
\begin{align*}
&\alpha_{i,j}=-q_{i,j}x_ix_j+\text{linear terms}\\
&\beta_{i}=x_i\sum_{k=1}^nq_{i,k}x_k+\text{linear terms}\\
&\gamma_{i}=x_i\sum_{k=1}^nq_{k,i}x_k+\text{linear terms}.
\end{align*}
Comparing the coefficient of $x_i^3x_j^3$ on both sides of \eqref{2.12a} 
we obtain 
$$q_{i,j}(q_{j,i}^2+a^2)=q_{j,i}(q_{i,j}^2+a^2),$$
which combined with $q_{i,j}+q_{j,i}=2a$ gives $q_{i,j}=q_{j,i}=a$. 
\end{proof}

In the rest of the paper, we will find the general solution of the 
compatibility conditions \eqref{2.12}, by equating the coefficients of 
the different powers of $x$, and we will determine the possible sets 
$V$, weights $W$ as well as explicit bases of orthogonal polynomials.

Notice that \eqref{3.4} implies that $\lambda_k$ must be at most quadratic 
in $k$. Thus, we have essentially two possible cases:
\begin{itemize}
\item[(i)] $\lambda_k$ is quadratic in $k$ (i.e. $a\neq 0$);
\item[(ii)] $\lambda_k$ is linear in $k$ (i.e. $a=0$).
\end{itemize}
Since we can always divide the equation $Du=\lambda_ku$ by a nonzero 
constant, we can take $a=-1$ in (i) and similarly we can normalize so that
$b=-1$ in (ii). These two cases are discussed in Sections \ref{se5} and 
\ref{se6} respectively. In the next section we treat in details the case 
$d=1$. 

\section{One dimensional orthogonal polynomials} \label{se4}

The one dimensional case is the simplest case with no compatibility 
conditions. Our goal in this section is to provide a detailed study in this 
simplest situation and to show how the procedure works (i.e. how we determine 
the weight function $W(x)$ and the set $V$). The discrete orthogonal 
polynomials of one variable will also serve as building blocks in the 
higher-dimensional cases. 

Historically, orthogonal polynomials satisfying difference equations of one 
variable were studied early in last century, we refer to \cite{Al,NSU} for 
references. It is well known that there are four families of such orthogonal 
polynomials of one variable, namely, Hahn polynomials, Meixner polynomials,
Krawtchouk polynomials, and Charlier polynomials. However, to our great
surprise, another family of orthogonal polynomials of discrete variable shows
up in our study.

For simplicity of notation, we will drop all indices within this section, i.e. 
$x=(x_1)\in\R$, $E=E_1$, $\fs=\fs_1$ etc. The operator $D$ will be 
of the form
\begin{equation}\label{4.1}
D=A(x)\fs\bs +B(x)\fs = \beta(x)E-(\beta(x)+A(x))\Id+A(x)E^{-1},
\end{equation}
where 
\begin{equation}\label{4.2}
A(x)=x(ax+l),\quad B(x)=bx+s, \text{ and }\beta(x)=A(x)+B(x).
\end{equation}
As we explained at the end of the previous section, we have two possible 
cases: \\
(i) quadratic eigenvalue $\lambda_k$, we can fix $a=-1$\\
(ii) linear eigenvalue $\lambda_k$, we can fix $a=0$, $b=-1$.


\subsection{Quadratic eigenvalue: $a=-1$} \label{sse4.1}
In order to make the 
formulas more symmetric, let us write $A(x)$,  $\beta(x)$ and $B(x)$ as 
\begin{subequations}\label{4.3}
\begin{align}
&A(x)=-x(x+\alpha_3)\label{4.3a}\\
&\beta(x)=-(x+\alpha_1+1)(x+\alpha_2+1)\label{4.3b}\\
&B(x)=\beta(x)-A(x)=(\alpha_3-\alpha_1-\alpha_2-2)x
-(1+\alpha_1)(1+\alpha_2),\label{4.3c}
\end{align}
\end{subequations}
i.e. we have replaced $l=-\alpha_3$, $b=\alpha_3-\alpha_1-\alpha_2-2$ and 
$s=-(1+\alpha_1)(1+\alpha_2)$. Thus \eqref{2.9} gives
\begin{equation}\label{4.4}
\frac{W(x)}{W(x-1)}=\frac{\beta(x-1)}{A(x)}
=\frac{(x+\alpha_1)(x+\alpha_2)}{x(x+\alpha_3)}.
\end{equation}

There are essentially 2 possibilities, depending on if the numerator in 
\eqref{4.4} vanishes for some positive integer $x$ or not.

\begin{Ha} Let \eqref{4.4} vanish for $x=N+1$ but be positive for 
$x\leq N$. Then one of $\alpha_1$ and $\alpha_2$ must be equal to $-N-1$, 
say $\alpha_2=-N-1$ and we can rewrite \eqref{4.4} as 
\begin{equation}\label{4.5}
\frac{W(x)}{W(x-1)}=\frac{(x+\alpha_1)(N+1-x)}{x(\beta_1+1+N-x)},
\end{equation}
where we denoted $\beta_1=-\alpha_3-N-1$. Thus, up to a constant factor, the 
weight becomes 
\begin{equation*}
W(x)=\binom{\alpha_1+x}{x}\binom{\beta_1+N-x}{N-x}.
\end{equation*}
The corresponding orthogonal polynomials on $V^1_N=\{0,1,\dots,N\}$ are the 
Hahn polynomials $Q_n(x;\alpha_1,\beta_1,N)$ given by 
\begin{equation}\label{4.6}
 Q_n(x;\alpha_1,\beta_1,N) = {}_3 F_2
\left(\begin{matrix} -n, n+\alpha_1+\beta_1+1, -x\\
       \alpha_1+1, -N \end{matrix}; 1\right). 
\end{equation}
In order to have a positive weight at all $x\in V^1_N$ the right-hand side of 
\eqref{4.5}
must be positive for $x=1,2,\dots,N$. Putting $x=1$ and $x=N$ we see that 
$(\alpha_1+1)(\beta_1+N)>0$ and $(\alpha_1+N)(\beta_1+1)>0$. Starting with 
$\beta_1> -1$ or $\beta_1 < -1$ shows that  there are two solutions,
\begin{itemize}
\item[(i)]  $\alpha_1>-1$,  $\beta_1>-1$, 
\item[(ii)] $\alpha_1<-N$, $\beta_1<-N$.
\end{itemize}
Conversely, it is clear that if $\alpha_1$ and $\beta_1$ satisfy (i) or 
(ii) the weight function $W(x)$ will have a constant sign on $V^1_N$. Hahn 
polynomials satisfy the following orthogonal relation
\begin{align}\label{4.7}
 \sum_{x=0}^N \frac{(\alpha_1+1)_x (\beta_1+1)_{N-x}}{x! (N-x)!}
   & Q_n(x;\alpha_1, \beta_1, N)Q_m(x;\alpha_1,\beta_1, N) \\
 = &\, \frac{(-1)^n n! (\beta_1+1)_n (n+\alpha_1+\beta_1+1)_{N+1}}
     {N! (2n+\alpha_1+\beta_1+1) (-N)_n (\alpha_1+1)_n} \delta_{n,m}. \notag
\end{align}
The difference operator $D$ takes the form
\begin{equation}\label{4.8}
D=x(\beta_1+N+1-x)\fs\bs+((\alpha_1+1)N-x(\alpha_1+\beta_1+2))\fs\\ .
\end{equation}
\end{Ha}

\begin{HaT} Assume now that the right hand side of \eqref{4.4} 
does not vanish for any $x\in\N$. Then it must be positive for all 
$x\in\N$. In order to have a positive function $W(x)$ and to have convergent 
series $\sum_{x\in\N_0}x^nW(x)$ for some $n\in\N_0$, the parameter 
$\alpha_3$ must be positive and  $\alpha_1,\alpha_2$ must satisfy one of the 
following 
\begin{itemize}
\item[(i)] $\alpha_1$ and $\alpha_2$ are both positive, or  there exists 
a negative integer $\kappa$ such that $\alpha_1,\alpha_2\in(\kappa,\kappa+1)$;
\item[(ii)] $\alpha_2=\bar{\alpha}_1\in\C\setminus\R$. 
\end{itemize}

For such $\alpha_1,\alpha_2,\alpha_3$ we can write the weight $W(x)$ as 
\begin{equation}\label{4.9}
W(x)=\frac{\binom{x+\alpha_1}{x}\binom{x+\alpha_2}{x}}
{\binom{x+\alpha_3}{x}}
=\frac{(\alpha_1+1)_x(\alpha_2+1)_x}{x!\,(\alpha_3+1)_x}.
\end{equation}
The series $\sum_{x\geq 0} x^nW(x)$ converges absolutely for 
$n<\alpha_3-\alpha_1-\alpha_2-1$, i.e. the corresponding polynomials 
will be orthogonal on $\N_0$ up to a given degree. 
It seems that these polynomials have not appeared in the literature before. 
Therefore, we  derive below their basic properties: explicit formula 
in terms of hypergeometric functions, the three term recurrence formula, 
etc.
\end{HaT}

\begin{Proposition}\label{pr4.1} 
The polynomials
\begin{equation}\label{4.10}
R_n(x)={}_3F_2
\left(\begin{matrix} -n , n-\alpha_3+\alpha_1+\alpha_2+1, -x \\
\alpha_1+1, \alpha_2+1 \end{matrix}\,; 1\right),
\end{equation}
satisfy
\begin{equation}\label{4.11}
D(R_n(x))=n(\alpha_3-\alpha_1-\alpha_2-1-n)R_n(x)
\end{equation}
where
\begin{equation}\label{4.12}
D=-x(x+\alpha_3)\bs\fs +((\alpha_3-\alpha_1-\alpha_2-2)x
-(1+\alpha_1)(1+\alpha_2))\fs,
\end{equation}
and the orthogonal relation
\begin{equation*}
\sum_{x=0}^{\infty}R_n(x)R_m(x)
\frac{(\alpha_1+1)_x(\alpha_2+1)_x}{x!\,(\alpha_3+1)_x}=0,
\end{equation*}
for $n\neq m$ such that $n+m<\alpha_3-\alpha_1-\alpha_2-1$.
\end{Proposition}

\begin{proof}
We look for a solution of \eqref{4.11} of the form
\begin{equation}\label{4.13}
R_n(x)=
\sum_{k=0}^na_k m_k(x),\quad\text{ where }\quad m_k(x)=\frac{(-x)_k}{k!}.
\end{equation}
Notice that 
\begin{equation}\label{4.14}
\fs m_k=-m_{k-1},\quad x\bs m_{k}=km_k, \text{ and }
xm_{k-1}=(k-1)m_{k-1}-km_k.
\end{equation}
Using these formulas we can easily calculate $D(R_n)$:
\begin{align*}
D(R_n(x))=& \sum_{k=1}^n [(x+\alpha_3)x\bs (a_km_{k-1}(x))
-((\alpha_3-\alpha_1-\alpha_2-2)x\\
&\qquad\qquad-(1+\alpha_1)(1+\alpha_2))a_km_{k-1}(x)]\\
=& \sum_{k=1}^n\left[(k-\alpha_3+\alpha_1+\alpha_2+1)x
+\alpha_3(k-1)+(1+\alpha_1)(1+\alpha_2)\right]a_km_{k-1}(x)\\
=& \sum_{k=0}^n\left[k(\alpha_3-\alpha_1-\alpha_2-1-k)a_k
+(\alpha_1+k+1)(\alpha_2+k+1)a_{k+1}\right],
\end{align*}
where in the last formula we assumed that $a_{n+1}=0$. Thus the difference 
equation \eqref{4.11} leads to the recursive relation
\begin{equation*}
a_{k+1}=-\frac{(n-k)(n+k-\alpha_3+\alpha_1+\alpha_2+1)}
{(\alpha_1+k+1)(\alpha_2+k+1)}a_{k},
\end{equation*}
whose solution with $a_0=1$ is 
\begin{equation}\label{4.15}
a_{k}=\frac{(-n)_k(n-\alpha_3+\alpha_1+\alpha_2+1)_k}
{(\alpha_1+1)_k(\alpha_2+1)_k}.
\end{equation}
This shows that the polynomials $R_n(x)$ defined by \eqref{4.10} satisfy 
\eqref{4.11}. The orthogonality follows from \reref{re3.4}.
\end{proof}

From the explicit formula we can derive other properties of these orthogonal 
polynomials.

\begin{Proposition}\label{pr4.2}
The polynomials $R_n$ satisfy the three-term relation
\begin{equation}\label{4.16}
xR_n(x)=\fa_{n}R_{n+1}(x)-(\fa_n+\fc_n)R_n(x)+\fc_nR_{n-1}(x),
\end{equation}
where the coefficients are given by
\begin{subequations}\label{4.17}
\begin{align}
&\fa_n=\frac{(n+\alpha_1+1)(n+\alpha_2+1)(n-\alpha_3+\alpha_1+\alpha_2+1)}
{(2n-\alpha_3+\alpha_1+\alpha_2+1)(2n-\alpha_3+\alpha_1+\alpha_2+2)}
                                                            \label{4.17a}\\
\intertext{and}
&\fc_n=-\frac{n(n-\alpha_3+\alpha_1)(n-\alpha_3+\alpha_2)}
{(2n-\alpha_3+\alpha_1+\alpha_2)(2n-\alpha_3+\alpha_1+\alpha_2+1)}.
                             \label{4.17b}
\end{align}
\end{subequations}
\end{Proposition}

\begin{proof}
Using again \eqref{4.14} we see that
\begin{equation*}
xR_n =\sum_{k=1}^nk(a_k-a_{k-1})m_k(x)-(n+1)a_{n}m_{n+1}(x),
\end{equation*}
where $a_k$ are given by \eqref{4.15}. It is now a simple matter to compare it 
with $\fa_{n}R_{n+1}(x)-(\fa_n+\fc_n)R_n(x)+\fc_nR_{n-1}(x)$ and verify that 
the coefficients of $m_k(x)$ agree for $k=0,1,\dots,n+1$.
\end{proof}

\begin{Proposition}\label{pr4.3}
The norm of the polynomial $R_n$ can be computed 
from the formula 
\begin{equation}\label{4.18}
\begin{split}
\langle R_n, R_n\rangle = &
\frac{n!}{(\alpha_3-\alpha_1-\alpha_2-2n-1)(1+\alpha_1)_n(1+\alpha_2)_n}\\
&\times \frac{\Gamma(\alpha_3-\alpha_1-\alpha_2-n)\Gamma(\alpha_3+1)}
{\Gamma(\alpha_3-\alpha_1-n)\Gamma(\alpha_3-\alpha_2-n)}.
\end{split}
\end{equation}
\end{Proposition}

\begin{proof} Using the three-term relation, a standard computation shows 
that
\begin{equation*}
\langle R_n, R_n\rangle =
\frac{c_nc_{n-1}\dots c_1}{a_{n-1}a_{n-2}\dots a_0}\langle R_0, R_0\rangle,
\end{equation*}
The quantity $\langle R_0, R_0\rangle=\langle 1, 1\rangle$ is 
a ${}_2F_1$ evaluated at $1$, which can be computed in closed form using 
Gauss' formula
\begin{equation*}
\langle 1, 1\rangle =
\frac{\Gamma(\alpha_3+1)\Gamma(\alpha_3-\alpha_1-\alpha_2-1)}
{\Gamma(\alpha_3-\alpha_1)\Gamma(\alpha_3-\alpha_2)}.
\end{equation*}
The statement follows from the above relations, \eqref{4.17} and the formula
$$
\Gamma(x)=(-1)^n(-x+1)_n\Gamma(x-n).
$$
\end{proof}

\begin{Remark} \label{re4.4}
Notice that although $\langle R_n,R_m\rangle$ converges only 
for finitely many $(n,m)$, formula \eqref{4.10} defines polynomials for 
every $n\in\N_0$. The operator $D$ given in \eqref{4.12} is admissible 
according to \deref{de3.3}. Moreover, equations \eqref{4.11} and 
\eqref{4.16} hold for every $n\in\N_0$. In particular, this implies that 
the polynomials $R_n(x)$ provide a solution to the discrete-discrete version 
of the bispectral problem discussed in \cite{Gr}.
\end{Remark}


\subsection{Linear eigenvalue: $a=0$ and $b=-1$}\label{sse4.2}
If $a=0$ and $b=-1$ equations \eqref{4.2} give $A(x)=lx,\quad B(x)=-x+s$, and 
$\beta(x)=A(x)+B(x)=(l-1)x+s$.
Thus
\begin{equation}\label{4.19}
\frac{W(x)}{W(x-1)}=\frac{(l-1)(x-1)+s}{lx}.
\end{equation}
We have now 2 possibilities: $l=1$ or $l\neq 1$.

\begin{Ch}
If $l=1$ equation \eqref{4.19} essentially 
reduces to 
\begin{equation*}
\frac{W(x)}{W(x-1)}=\frac{s}{x},
\end{equation*}
i.e. $s$ must be a positive number and (up to a constant factor) 
$W(x)=s^x/x!$. The corresponding polynomials are the Charlier polynomials 
$$C_n(x;s)={}_2F_0
\left(\begin{matrix} -n , -x \\
-\end{matrix}\,; -\frac{1}{s}\right)$$ 
orthogonal on $\N_0$. The difference operator takes the form
\begin{equation}\label{4.20}
D=x\fs\bs +(s-x)\fs=sE-(x+s)\Id+xE^{-1}.
\end{equation}
\end{Ch}

Assume now that $l\neq 1$. Then, we can rewrite \eqref{4.19} as
\begin{equation}\label{4.21}
\frac{W(x)}{W(x-1)}=c\frac{x+d}{x},
\end{equation}
where we have put $c=(l-1)/l\neq 0$ and $d=(1+s-l)/(l-1)$.

We now have two possibilities depending on the sign of $c$.

\begin{Kr} 
If $c<0$ then \eqref{4.21} can be 
rewritten as 
\begin{equation}\label{4.22}
\frac{W(x)}{W(x-1)}=(-c)\frac{-d-x}{x},
\end{equation}
from which it follows that $(-d)$ must be a positive integer (otherwise, 
$W(x)$ will become negative at some point). If we put $d=-N-1$ and 
$p=c/(c-1)\in(0,1)$ we obtain the Krawtchouk polynomials 
$$K_n(x;p,N)={}_2F_1
\left(\begin{matrix} -n , -x \\
-N\end{matrix}\,; \frac{1}{p}\right)$$
orthogonal on $V^1_N=\{0,1,\dots,N\}$ with respect to the weight 
$W(x)=p^x(1-p)^{N-x}\binom{N}{x}$.
The difference operator takes the form 
\begin{equation}\label{4.23}
D=x(1-p)\fs\bs +(Np-x)\fs.
\end{equation}
\end{Kr}

\begin{Me}
If $c>0$ then $W(1)>0$ implies that 
$\beta:=d+1>0$. In order to have convergent series, we need $c<1$. This 
leads to the Meixner polynomials 
\begin{equation}\label{4.24}
M_n(x;\beta,c)={}_2F_1
\left(\begin{matrix} -n , -x \\
\beta\end{matrix}\,; 1-\frac{1}{c}\right),
\end{equation}
depending on two parameters $0<c<1$ and $\beta>0$, orthogonal on $\N_0$ 
with respect to the weight $W(x)=c^x(\beta)_x/x!$. The difference operator is
\begin{equation}\label{4.25}
D=\frac{x}{1-c}\fs\bs +\left(-x-\frac{c\beta}{1-c}\right)\fs.
\end{equation}
The Meixner polynomials satisfy the orthogonal relation
\begin{equation} \label{4.26}
  \sum_{x=0}^\infty  \frac{(\beta)_x}{x!}c^x M_n(x,\beta,c)M_m(x,\beta,c)
      = \frac{c^{-n} n!}{(\beta)_n(1-c)^\beta} \delta_{m,n}.
\end{equation}
\end{Me}


\section{Multivariable case with quadratic eigenvalue: $a=-1$}\label{se5}

We first solve the compatibility conditions \eqref{2.12} in arbitrary 
dimension $d>1$, and then we write explicitly all possible 
orthogonal polynomials with the corresponding weights.
For $a=-1$ formulas \eqref{3.16}-\eqref{3.18} give:
\begin{subequations}\label{5.1}
\begin{align}
&\alpha_{i,j}=x_j(x_i-l_{i,j})\label{5.1a}\\
&\beta_i=-x_i\sum_{k=1}^dx_k+\sum_{k=1}^dl_{i,k}x_k+bx_i+s_i\label{5.1b}\\
&\gamma_i=x_i\left(-\sum_{k=1}^dx_k+\sum_{k=1}^dl_{k,i}\right).\label{5.1c}
\end{align}
\end{subequations}
The next proposition gives the general solution to the compatibility 
conditions \eqref{2.12}
\begin{Proposition} For $a=-1$ there are 2 solutions to the compatibility 
conditions. \\

Solution 1:
\begin{subequations}\label{5.2}
\begin{align}
&\alpha_{i,j}=x_j(x_i-l_{i})\label{5.2a}\\
&\beta_i=-(x_i-l_i)\left(\sum_{k=1}^dx_k-b-r\right)\label{5.2b}\\
&\gamma_i=-x_i\left(\sum_{k=1}^dx_k-\sum_{k=1}^dl_{k}-r\right),\label{5.2c}
\end{align}
\end{subequations}
where $b$, $r$, and $\{l_i\}_{i=1}^d$ are free parameters.\\

Solution 2:
\begin{subequations}\label{5.3}
\begin{align}
&\alpha_{i,j}=x_j(x_i-l_{i})\\
&\beta_i=-(x_i-l_i)\left(\sum_{k=1}^dx_k-\sum_{k=1}^dl_k+1-r_i\right)\\
&\gamma_i=-x_i\left(\sum_{k=1}^dx_k-\sum_{k=1}^dl_{k}-r_i\right),
\end{align}
\end{subequations}
where $\{l_i,r_i\}_{i=1}^d$ are free parameters and $b=\sum_{k=1}^dl_k-1$.\\
\end{Proposition}

\begin{proof}
Let us take $i\neq j$ and let $k\neq i,j$. 
Comparing the coefficients of $x_i^3x_jx_k$ and $x_i^3x_j$ on both sides of 
\eqref{2.12a} we get:
\begin{equation}\label{5.4}
-(l_{j,k}+l_{j,i}+2)=-2(l_{j,i}+1), 
\end{equation}
and
\begin{equation}\label{5.5}
\begin{split}
&(l_{j,i}+1)\left(\sum_{k=1}^d l_{k,i}\right)-(s_j-1-b-l_{j,j})+
(1+l_{i,j})(1+l_{j,i})=\\
&\qquad -(l_{j,i}+1)
\left[-\left(\sum_{k=1}^d l_{k,j}\right)-(b+l_{i,i}+2)\right].
\end{split}
\end{equation}
equation \eqref{5.4} simply means that for $i\neq j$, $l_{j,i}$ is 
independent of $i$, i.e. we can put $l_{j}:=l_{j,i}$ for $j\neq i$. 
Using this in \eqref{5.5} we obtain that 
\begin{equation}
s_j=l_j(l_j-l_{j,j}-b).
\end{equation}

Our formulas simplify as follows
\begin{subequations}\label{5.7}
\begin{align}
&\alpha_{i,j}=x_j(x_i-l_{i})\\
&\beta_i=-x_i\sum_{k=1}^dx_k+l_i\sum_{k\neq i} x_k+(l_{i,i}+b)x_i
-l_i(l_{i,i}+b-l_i)\nonumber\\
&\quad =-(x_i-l_i)\left(\sum_{k=1}^dx_k-l_{i,i}-b+l_i\right)\\
&\gamma_i=x_i\left(-\sum_{k=1}^dx_k+\sum_{k\neq i}l_{k}+l_{i,i}\right).
\end{align}
\end{subequations}

With the above formulas it is very easy to solve completely \eqref{2.12a}. 
Indeed, after canceling the common factor $x_ix_j(x_i-l_i-1)(x_j-l_j-1)$ 
equation \eqref{2.12a} gives
\begin{equation}\label{5.8}
(l_{i,i}-l_{i}-l_{j,j}+l_{j})\left(\sum_{k=1}^d l_k-b-1\right)=0.
\end{equation}
Let us denote $r_i=l_{i,i}-l_{i}$. Then the first factor in the 
last formula is simply $r_i-r_j$.

From \eqref{5.8} it follows that we have 2 possibilities: $r_i=r_j$ for 
all $i\neq j$, i.e. $r:=r_i$ is independent of $i$, or 
$b=\sum_{k=1}^dl_k-1$. The first leads to Solution 1, and the second leads 
to Solution 2. It is immediate to see that the functions defined by 
\eqref{5.2} or \eqref{5.3} satisfy \eqref{2.12b}.
\end{proof}

Next, we want to determine the possible weights and orthogonal polynomials.
It is easy to see that there is no weight corresponding to Solution 2. 
Otherwise, from \eqref{5.3} and \eqref{2.9} we will have
\begin{equation}\label{5.9}
\frac{W(x)}{W(x-e_i)}=\frac{x_i-l_i-1}{x_i}\text{ for }i=1,2,\dots,d.
\end{equation}
Then, for $x=e_i$ it follows that $l_i<0$, i.e we can put $-l_i=r_i>0$ and 
therefore, up to a constant factor, $W(x)$ is given by
\begin{equation*}
W(x)=\prod_{i=1}^d\frac{(r_i)_{x_i}}{x_i!},
\end{equation*}
but then even $\langle 1,1\rangle $ will diverge.  

Let us now concentrate on Solution 1, given by \eqref{5.2}. 
Equations \eqref{5.2} and \eqref{2.9} give
\begin{equation}\label{5.10}
\frac{W(x)}{W(x-e_i)}=\frac{x_i-l_i-1}{x_i}
\frac{\left(\sum_{k=1}^dx_k-b-r-1\right)}
{\left(\sum_{k=1}^dx_k-\sum_{k=1}^dl_{k}-r\right)},
\end{equation}
for $i=1,2,\dots,d$.
For a vector $y\in\R^k$ we will denote 
\begin{equation}\label{5.11}
|y|=\sum_{i=1}^ky_i.
\end{equation}
If we put 
\begin{equation}\label{5.12}
\sigma_i=-(l_i+1),
\end{equation} 
and 
\begin{equation}\label{5.13}
\sigma =(\sigma_1,\sigma_2,\dots,\sigma_d)
\end{equation}
we can write \eqref{5.10} as
\begin{equation}\label{5.14}
\frac{W(x)}{W(x-e_i)}=\frac{x_i+\sigma_i}{x_i}
\frac{\left(|x|-b-r-1\right)}
{\left(|x|+|\sigma|+d-r\right)}.
\end{equation}
We have several possible sub-cases, depending on whether the right hand side  
in equation \eqref{5.14} can vanish or not. We divide them first in 
two big subclasses:
\begin{itemize}
\item $b+r\notin\N_0$;
\item $b+r\in\N_0$.
\end{itemize}


\subsection{Case 1: $b+r\not\in\N_0$}\label{sse5.1}
Then $|x|-b-r-1$ does not vanish. We have two possibilities:


\subsubsection{} For some $i\in\{1,2\dots,d\}$, $(-\sigma_i)\not\in\N$. 
We show below that this implies that $\sigma_j>-1$ for all 
$j\in\{1,2\dots,d\}$. Take $j\neq i$, arbitrary $k\in\N$ and $n\in\N_0$.
From \eqref{5.14} we see that 
\begin{equation}\label{5.15}
\frac{W(ke_i+ne_j)}{W((k-1)e_i+ne_j)}=\frac{k+\sigma_i}{k}
\frac{k+n-b-r-1}{k+n+|\sigma|+d-r}>0.
\end{equation}
The second ratio is positive for large $n$, and therefore $k+\sigma_i>0$ 
for every $k\in\N$, i.e. $\sigma_i>-1$. Equation \eqref{5.15} for $n=0$ 
implies that
$$\frac{k-b-r-1}{k+|\sigma|+d-r}>0,$$
for every $k\in\N$. But then for every $j\in\{1,2,\dots,d\}$ we must have
\begin{equation*}
\frac{W(e_j)}{W(0)}=\frac{1+\sigma_j}{1}
\frac{1-b-r-1}{1+|\sigma|+d-r}>0,
\end{equation*}
showing that $\sigma_j>-1$.

Since $\sigma_j>-1$ for all $j$ and $b+r\notin\N_0$ we must have $V=\N_0^d$.
Up to a constant, the weight $W(x)$ is given by
\begin{equation*}
W(x)=\prod_{i=1}^d\frac{(\sigma_i+1)_{x_i}}{x_i!}\frac{(-b-r)_{|x|}}
{(|\sigma|+d-r+1)_{|x|}}
\end{equation*}
Changing the parameters
\begin{equation*}
\gamma=|\sigma|+d-r\qquad \beta=-(b+r)-1,
\end{equation*}
the weight takes the form 
\begin{equation*}
W(x)=\prod_{i=1}^d\frac{(\sigma_i+1)_{x_i}}{x_i!}\frac{(\beta+1)_{|x|}}
{(\gamma+1)_{|x|}},
\end{equation*}
and the free parameters are $\{\sigma_i\}$, $\beta$ and $\gamma$.
The difference operator becomes
\begin{equation}\label{5.16}
\begin{split}
D&=-\sum_{1\leq i,j\leq d}x_j[\sigma_i+1+x_i+\delta_{i,j}(\gamma-|\sigma|-d)]
\fs_i\bs_j\\
&\qquad+
\sum_{i=1}^d[x_i(\gamma-\beta-|\sigma|-d-1)-(1+\beta)(1+\sigma_i)]\fs_i,
\end{split}
\end{equation}
and $\lambda_n=n(-n+\gamma-\beta-|\sigma|-d)$.

For a vector $y=(y_1,y_2,\dots,y_d)\in\R^d$ we will denote
\begin{equation}\label{5.17}
y^j=(y_j,y_{j+1},\dots,y_d)\text{ and }Y_j=(y_1,y_2,\dots,y_j),
\end{equation}
with the convention that $y^{d+1}=0$ and $Y_0=0$. 
Let us denote by $R_n(x;\alpha_1,\alpha_2,\alpha_3)$ the orthogonal 
polynomials given by \eqref{4.10}.

\begin{Theorem} \label{th5.2}
For $\nu \in \N_0^d$, such that
$2|\nu| < \gamma - |\sigma|-\beta - d -1$, the 
polynomials 
\begin{equation*}
R_\nu(x;\sigma,\beta,\gamma) = 
    \prod_{j=1}^d  (\alpha_{2,j}+1)_{\nu_j}
  R_{\nu_j}(x_j; \sigma_j,\alpha_{2,j},\alpha_{3,j}),        
\end{equation*}
where $\alpha_{2,j}$ and $\alpha_{3,j}$ are given by
\begin{align*}
&\alpha_{2,j}=\beta+|\nu^{j+1}|+|X_{j-1}|\\
&\alpha_{3,j}=\gamma-|\nu^{j+1}|-|\sigma^{j+1}|-(d-j)+|X_{j-1}|
\end{align*}
satisfy the difference equation
$$D\phi_{\nu}=\lambda_{|\nu|}\phi_{\nu}$$ 
and the orthogonal relation
\begin{equation} \label{5.18}
\sum_{x \in \N_0^d}  R_\nu(x; \sigma,\beta,\gamma) 
                    R_\mu(x; \sigma,\beta,\gamma)
\prod_{i=1}^d \frac{(\sigma_i+1)_{x_i}}{x_i!} 
   \frac{(\beta+1)_{|x|}}{(\gamma+1)_{|x|}} 
 = A_\nu \delta_{\nu,\mu},
\end{equation}
where 
\begin{equation}\label{5.19}
\begin{split}
A_\nu=&\frac{(1+\beta)_{|\nu|}\Gamma(1+\gamma)
\Gamma(\gamma-\beta-|\sigma|-2|\nu|-d)}
{\Gamma(\gamma-|\sigma|-|\nu|+1-d)\Gamma(\gamma-\beta)}\\
&\times
 \prod_{j=1}^d \frac{\nu_j! (\alpha_{3,j}-\alpha_{2,j} - \nu_j)_{\nu_j}
     (\alpha_{3,j-1}-\alpha_{2,j-1} +1)_{\nu_j}} {(1+\sigma_j)_{\nu_j}}.
\end{split}
\end{equation}
\end{Theorem}

\begin{proof}
Notice that $\alpha_{3,j}-\alpha_{2,j}$ is independent of $x$, 
and therefore from \eqref{4.10} one can immediately see that $\phi_{\nu}$ is 
indeed a polynomial of $x$ of total degree $|\nu|$. From \thref{th3.7} we 
know that for every $k\in\N_0$, the equation $Du=\lambda_ku$ has 
$r_k=\binom{k+d-1}{k}=\dim(\Pi^d_{k}/\Pi^d_{k-1})$ linearly independent 
solutions. 
Therefore, it is enough to prove that \eqref{5.18} holds.

In \eqref{5.18} we will first sum with respect to $x_d$, then with respect to 
$x_{d-1}$, and so on. Writing
$$
(\beta+1)_{|x|} = (\beta+1+|X_{d-1}|)_{x_d} (\beta+1)_{|X_{d-1}|}  
$$
and extracting only the terms depending on $x_d$ we get the sum
\begin{align*}
I_{\nu_d,\mu_d}: &= \sum_{x_d \ge 0} R_{\nu_d}(x_d)R_{\mu_d}(x_d) 
    \frac{(\sigma_d+1)_{x_d}}{x_d!} \frac{(\beta+1+|X_{d-1}|)_{x_d}}
     {(\gamma+1+|X_{d-1}|)_{x_d}}\\
&=\sum_{x_d \ge 0} R_{\nu_d}(x_d)R_{\mu_d}(x_d) 
    \frac{(\sigma_d+1)_{x_d}}{x_d!} \frac{(\alpha_{2,d}+1)_{x_d}}
     {(\alpha_{3,d}+1)_{x_d}},
\end{align*}
where $R_m(x_d)=R_{m}(x_d;\sigma_d,\alpha_{2,d},\alpha_{3,d})$.
Using now \eqref{4.18} and the fact that $\alpha_{3,d}-\alpha_{2,d}=
\gamma-\beta$ is independent of $x$, we see that 
$$
 I_{\nu_d,\mu_d} 
= \delta_{\nu_d,\mu_d} B_{\nu_d}
\frac{\Gamma(\gamma+|X_{d-1}|+1)}  
{ (1+\beta +|X_{d-1}|)_{\nu_d}\Gamma(\gamma-\sigma_d-\nu_d+|X_{d-1}|)},
$$
where $B_{\nu_d}$ is a constant (independent of $x$), whose value 
can be extracted from \eqref{4.18}
\begin{equation*}
B_{\nu_d}=\frac{\nu_d!\Gamma(\gamma-\beta-\sigma_d-\nu_d)}
{(1+\sigma_d)_{\nu_d}(\gamma-\beta-\sigma_d-2\nu_d-1)
\Gamma(\gamma-\beta-\nu_d)}.
\end{equation*}
Using the fact that 
\begin{equation}\label{5.20}
\Gamma(x)=(x-n)_n\Gamma(x-n),
\end{equation}
we see that
$$
 I_{\nu_d,\mu_d} 
= \delta_{\nu_d,\mu_d} C_{\nu_d}
\frac{(\gamma+1)_{|X_{d-1}|}}  
{ (1+\beta +|X_{d-1}|)_{\nu_d}(\gamma-\sigma_d-\nu_d)_{|X_{d-1}|}},
$$
with $C_{\nu_d}=B_{\nu_d} {\Gamma(\gamma+1)}/{\Gamma(\gamma-\sigma_d-\nu_d)}$.
This shows, in particular, that
\begin{equation}\label{5.21}
\begin{split}
& [ (\beta+1+|X_{d-1}|)_{\nu_d}]^2 
\frac{(\beta+1)_{|X_{d-1}|}}{(\gamma+1)_{|X_{d-1}|}} I_{\nu_d,\mu_d}\\
 &  \qquad\qquad = \delta_{\nu_d,\mu_d}    C_{\nu_d}  (\beta+1)_{\nu_d} 
   \frac{(\beta+1+\nu_d)_{|X_{d-1}|}} {(\gamma -\sigma_d-\nu_d)_{|X_{d-1}|}},
\end{split}
\end{equation}
where we have used the identity 
\begin{equation}\label{5.22}
(\beta+1)_{|X_{d-1}|} (\beta+1+ |X_{d-1}|)_{\nu_d}  = (\beta+1)_{|\nu_d|} 
(\beta+1+\nu_d)_{|X_{d-1}|}. 
\end{equation}
Equation \eqref{5.21} shows that the remaining $d-1$ fold sums of 
$\langle R_\nu,R_\mu\rangle$ have exactly the same structure as 
that of the original $d$ fold sums with $\beta$ and $\gamma$ replaced by 
$\beta+\nu_d$ and $\gamma-\sigma_d-\nu_d-1$, respectively.  
In other words, it shows that we can use induction to complete the proof. 
For $A_{\nu}$ we obtain the following formula:
\begin{equation}\label{5.23}
\begin{split}
A_\nu=&\frac{(1+\beta)_{|\nu|}\Gamma(1+\gamma)}
{\Gamma(\gamma-|\sigma|-|\nu|+1-d)}
\prod_{j=1}^d \frac{\nu_j!}{(1+\sigma_j)_{\nu_j}
(\gamma-\beta-|\sigma^j|-2|\nu^j|+j-d-1)}\\
&\quad \times 
\prod_{j=1}^d \frac{\Gamma(\gamma-\beta-|\sigma^j|-2|\nu^{j+1}|-\nu_j+j-d)}
{\Gamma(\gamma-\beta-|\sigma^{j+1}|-2|\nu^{j+1}|-\nu_j+j-d)}.
\end{split}
\end{equation}
Using several times \eqref{5.20} one can rewrite the right-hand side as 
in \eqref{5.19}.
\end{proof}


\subsubsection{} For every $i\in\{1,2\dots,d\}$ we have 
$(-\sigma_i)=l_i+1 \in\N$, i.e. $l_i\in\N_0$. The corresponding 
$V$ is the parallelepiped
\begin{equation*}
V^d_l=\{x\in\N_0^d: x_i\leq l_i\}.
\end{equation*}
By \eqref{5.10} the weight function in this case is given by 
\begin{equation} \label{5.24}
W(x) = \prod_{i=1}^d 
\frac{(- l_i)_{x_i} }{x_i!} \frac{(\beta+1)_{|x|}}{(-|l |-r+1)_{|x|}}, 
\end{equation}
where we again set $\beta = - (b+r) -1$. Putting $x=e_i$ and 
$x=(l_1,l_2,\dots,l_d)$ in \eqref{5.10} we get that
$$\frac{\beta+1}{|l|+r-1}>0
\quad\text{ and }\quad
\frac{|l|+\beta}{r}>0.$$
From this we see that the parameters $\beta$ and $r$ must satisfy one of the 
following conditions:
\begin{itemize}
\item[(i)] $\beta>-1$ and $r>0$;
\item[(ii)] $\beta<-|l|$ and $r<-|l|+1$.
\end{itemize}
Recall that $b + r \notin \N_0$, but everything will hold even if 
$b+r=-\beta-1\in\N_0$ as long as $\beta+|l|<0$ (i.e. if (ii) holds),
 because $(\beta+1)_{|x|}$ will not vanish for $x\in V^d_l$.

In the following we will use the notations $y^j$ and $Y_j$ defined in 
\eqref{5.17}. For example, $|L_j| = l_1 + \cdots +l_j$.  If
$\nu, l \in \N_0^d$, then $\nu\le l$ means $\nu_i\le l_i$ for $1 \le i \le d$.
Recall that the Hahn polynomials are denoted by $Q_n(x;\alpha_1, \beta_1, N)$. 

\begin{Theorem}\label{th5.3}
Let $l_i\in \N_0$, $1 \le i \le d$. For $\nu \in \N_0^d$, $\nu_i\le l_i$, the 
polynomials
\begin{equation} \label{5.25}
\phi_\nu(x;\beta,r,l)= 
   \prod_{i=1}^d (\alpha_{1,j}+1)_{\nu_j} Q_{\nu_j}(x_j; \alpha_{1,j},   
            \alpha_{2,j},l_j),
\end{equation}
where $\alpha_{1,j}$ and $\alpha_{2,j}$ are given  by 
$$
   \alpha_{1,j}  = \beta+ |\nu^{j+1}|+|X_{j-1}| \quad \hbox{and}\quad 
   \alpha_{2,j}  =|L_{j-1}| - |X_{j-1}|+|\nu^{j+1}|+ r -1,
$$
satisfy the difference equation
$$
     D \phi_\nu = \lambda_{|\nu|}\phi_\nu
$$
and the orthogonal relation
\begin{equation}\label{5.26}
 \sum_{\nu\le l} \phi_{\nu}(x;\beta,r,l)\phi_{\mu}(x;\beta,r,l) \prod_{i=1}^d
   \frac{(- l_i)_{x_i} }{x_i!} \frac{(\beta+1)_{|x|}}{(-|l |-r+1)_{|x|}} 
  = B_\nu \delta_{\nu,\mu},
\end{equation}
where the normalization constant is given by
\begin{equation}\label{5.27}
B_\nu =  \frac{(-1)^{|\nu|}(1+\beta)_{|\nu|}}{(r+|\nu|)_{|l|-|\nu|}}
\prod_{j=1}^d\frac{\nu_j!(\beta+r+2|\nu^{j+1}|+\nu_j+|L_{j-1}|)_{l_j+1}}
{(-l_j)_{\nu_j}(\beta+r+2|\nu^{j}|+|L_{j-1}|)}.
\end{equation}
\end{Theorem}

\begin{proof}
Since $\alpha_{1,j} +\alpha_{2,j}$ is independent of $x$, it is easy to see 
from \eqref{4.6} that $\phi_\nu$ is indeed a polynomial of $x $ of total 
degree $|\nu|$. We proceed as in the proof of \thref{th5.2}. In \eqref{5.26} 
we will first sum with respect to $x_d$. Using the fact that 
\begin{align*}
      (-|l|-r+1)_{|x|} &= (-|l|-r+1+|X_{d-1}|)_{x_d}(-|l|-r+1)_{|X_{d-1}|}, \\
        (\beta+1)_{|x|}& = (\beta+1)_{|X_{d-1}|}(\beta+1+|X_{d-1}|)_{x_d},
\end{align*}
we can split the weight function as a product, 
$$  
 W(x)=\frac{(-l_d)_{x_d}(\beta+1+|X_{d-1}|)_{x_d}}
     {x_d!(-|l|-r+1+|X_{d-1}|)_{x_d}}W'(X_{d-1})
$$
where 
$$
W'( X_{d-1} ) = \prod_{i=1}^{d-1} \frac{(- l_i)_{x_i} }{x_i!} 
     \frac{(\beta+1)_{|X_{d-1}|}}{(-|l |-r+1)_{|X_{d-1}|}}. 
$$
It is easy to verify that 
$$
 \frac{ (-l_d)_{x_d} (\beta+1+|X_{d-1}|)_{x_d}}{x_d!(-|l|-r+1+|X_{d-1}|)_{x_d}}
     = \frac{l_d!}{(\alpha_{2,d}+1)_{l_d}}
      \binom{x_d+\alpha_{1,d}}{x_d}
          \binom{l_d- x_d+ \alpha_{2,d}}{l_d-x_d}. 
$$
Hence, using the fact that $\alpha_{1,d}$ and $\alpha_{2,d}$ are independent of
$x_d$, the sum over $x_d$ in \eqref{5.26} becomes 
\begin{align*}
I_{\nu_d,\mu_d}:=& \frac{l_d!}{(\alpha_{2,d}+1)_{l_d}}(\alpha_{1,d}+1)_{\nu_d}
        (\alpha_{1,d}+1)_{\mu_d} \\
    \times &\sum_{x_d=0}^{l_d}  Q_{\nu_d}(x_d) Q_{\mu_d}(x_d)
        \binom{x_d+\alpha_{1,d}}{x_d} \binom{l_d- x_d+ \alpha_{2,d}}{l_d-x_d},
\end{align*}
where $Q_m(x_d) = Q_m(x_d,\alpha_{1,d}, \alpha_{2,d}, l_d)$. Using 
\eqref{4.7} and simplifying, we get
$$
    I_{\nu_d,\mu_d} = \delta_{\nu_d,\mu_d} 
      \frac{(-1)^{\nu_d} \nu_d! (\nu_d+\beta+ |L_{d-1}|+r)_{l_d+1}}
           { (-l_d)_{\nu_d} (2\nu_d+\beta+|L_{d-1}|+r) } 
 \frac{(\alpha_{2,d}+1)_{\nu_d}}{(\alpha_{2,d}+1)_{l_d}}
     (\alpha_{1,d}+1)_{\nu_d}.$$  
In the above expression, the term $(\alpha_{1,d}+1)_{\nu_d}$ contains variables
$x_1,\ldots, x_{d-1}$. Combining this term with $W'(X_{d-1})$ and using 
\eqref{5.22}, we see that the weight function for $x_1, \ldots,x_{d-1}$ becomes
$$
W_{d-1}(X_{d-1})=(\beta+1)_{\nu_d}\prod_{i=1}^{d-1}\frac{(- l_i)_{x_i} }{x_i!} 
         \frac{(\beta+1+\nu_d)_{|X_{d-1}|}}{(-|l |-r+1)_{|X_{d-1}|}}
            \frac{(\alpha_{2,d}+1)_{\nu_d}}{(\alpha_{2,d}+1)_{l_d}}. 
$$  
By the definition of $\alpha_{2,d}$ and expanding the Pochhammer symbols 
gives 
\begin{align*}
&(-|l |-r+1)_{|X_{d-1}|} 
  \frac{(\alpha_{2,d}+1)_{l_d}}{(\alpha_{2,d}+1)_{\nu_d}} \\
 & \qquad   =(-1)^{l_d - \nu_d}  (-|l |-r+1) (-|l |-r+2)\ldots 
(-|L_{d-1} |-r-\nu_d+|X_{d-1}|) \\
 & \qquad 
      =  (-|L_{d-1}|-r+1-\nu_d)_{|X_{d-1}|}(|L_{d-1}|+r+\nu_d)_{l_d-\nu_d}.  
\end{align*}
Consequently, the weight function $W_{d-1}$ becomes 
$$
 W_{d-1}(X_{d-1}) = \frac{(\beta+1)_{\nu_d}}{(|L_{d-1}|+r+\nu_d)_{l_d-\nu_d}} 
     \prod_{i=1}^{d-1} \frac{(- l_i)_{x_i} }{x_i!} 
      \frac{(\beta+1+\nu_d)_{|X_{d-1}|}}{(-|L_{d-1} |-r+1-\nu_d)_{|X_{d-1}|}}.
$$
Apart from a constant multiple, the weight function $W_{d-1}$ has exactly the 
same structure of $W(x)$ with $\beta$ replaced by $\beta+\nu_d$, $r$ replaced 
by $r + \nu_d$, $l$ replaced by $L_{d-1}$, respectively, and one variable 
less. Hence, we can use induction to complete the proof.   
For $B_{\nu}$ we obtain the following formula
\begin{equation*}
B_\nu =  (-1)^{|\nu|}(1+\beta)_{|\nu|}
\prod_{j=1}^d\frac{\nu_j!(\beta+r+2|\nu^{j+1}|+\nu_j+|L_{j-1}|)_{l_j+1}}
{(-l_j)_{\nu_j}
    (\beta+r+2|\nu^{j}|+|L_{j-1}|)(|L_{j-1}|+r+|\nu^j|)_{l_j-\nu_j}},
\end{equation*}
which combined with 
$$\prod_{j=1}^d(|L_{j-1}|+r+|\nu^j|)_{l_j-\nu_j}=(r+|\nu|)_{|l|-|\nu|},$$
gives \eqref{5.27}.
\end{proof}


\subsection{Case 2: $b+r=N\in\N_0$} \label{sse5.2}
Again we have two possibilities: 


\subsubsection{}  For every $i\in\{1,2\dots,d\}$, 
$l_i+1\notin\{1,2,\dots,N\}$. In this case the numerator in \eqref{5.14} is 
zero when $|x| = N+1$ and we have $V = V^d_N$.  The corresponding orthogonal 
polynomials are the Hahn  polynomials studied in \cite{KM}.  Setting 
$\sigma_i = -l_i -1$ for $1\le i \le d$, 
$\sigma = (\sigma_1,\ldots,\sigma_{d+1})$,  
\begin{equation}\label{5.28}
  r = N+|\sigma|+d+1, \quad \hbox{and} \quad   b = -(|\sigma|+d+1),
\end{equation}
equation \eqref{5.14} gives
\begin{equation*}
\frac{W(x)}{W(x-e_i)}=\frac{x_i+\sigma_i}{x_i}
\frac{N+1-|x|}{N+1+\sigma_{d+1}-|x|}.
\end{equation*}
The above ratio must be positive for all $x\in V^d_N$. In particular, for 
$x=e_i$ and $x=Ne_i$ we see that
$$\frac{\sigma_i+1}{N+\sigma_{d+1}}>0\quad \text{ and }\quad
\frac{\sigma_i+N}{1+\sigma_{d+1}}>0 \quad\text{ for }i=1,2\dots,d.$$
From this it follows easily that the parameters $\{\sigma_{i}\}_{i=1}^{d+1}$ 
satisfy one of the following conditions:
\begin{itemize}
\item[(i)]  $\sigma_i>-1$ for $i=1,2,\dots,d+1$;
\item[(ii)] $\sigma_i<-N$ for $i=1,2,\dots,d+1$.
\end{itemize}

The weight function takes the form
\begin{equation} \label{5.29}
W(x) = \prod_{i=1}^d \binom{x_i+\sigma_i}{x_i} 
\binom{N-|x|+\sigma_{d+1}}{N-|x|}
\end{equation}
with $\{\sigma_i\}$ as free parameters. If $\sigma_i>-1$ for 
$i=1,2,\dots,d+1$ or if 
$\sigma_i<-N$ for $i=1,2,\dots,d+1$ but $N$ is even then 
$W(x)>0$ on $V^d_N$. If $\sigma_i<-N$ for $i=1,2,\dots,d+1$ and if $N$ is odd
then $W(x)<0$ on $V^d_N$, so one needs to change the sign in formula 
\eqref{5.29}, in order to get a positive function.

The difference operator takes the form 
\begin{align*}
&  D = \sum_{i=1}^d x_i (N-x_i+|\sigma|-\sigma_i+d)\fs_i\bs_i  -
     \sum_{1 \le i\ne j \le d}x_j(x_i+\sigma_i+1)\fs_i\bs_j \\
& \qquad\qquad 
+\sum_{i=1}^d [(N-x_i)(\sigma_i+1)- x_i(|\sigma|-\sigma_i+d)] \fs_i 
 \end{align*}
and the eigenvalues are $\lambda_n=  - n (n+|\sigma|+d)$. 

\begin{Theorem}\label{th5.4}
For $\nu \in \N_0^d$ and $|\nu|\le N$, the polynomials 
 \begin{align} \label{5.30}
 Q_\nu(x;\sigma, N) = &  \frac{(-1)^{|\nu|}}{(-N)_{|\nu|}}
  \prod_{j=1}^d  \frac{(\sigma_j+1)_{\nu_j}}{(a_j+1)_{\nu_j}}
    (-N + |X_{j-1}|+|\nu^{j+1}|)_{\nu_j} \\
 & \qquad\quad  \times Q_{\nu_j}(x_j; \sigma_j, a_j, N-|X_{j-1}|-|\nu^{j+1}|),
\notag
\end{align}
where $a_j = |\sigma^{j+1}|+2 |\nu^{j+1}| + d-j$, satisfy the difference 
equation
$$
     D Q_\nu = \lambda_{|\nu|}Q_\nu
$$
and the orthogonal relation
\begin{equation} \label{5.31}
\sum_{|x|\le N} Q_\nu(x; \sigma,N) Q_\mu(x; \sigma, N)
\prod_{i=1}^{d}\binom{x_i+\sigma_i}{x_i} \binom{N-|x|+\sigma_{d+1}}{N-|x|}
 = A_\nu \delta_{\nu,\mu}, 
\end{equation}
where $A_\nu$ is given by
\begin{equation}\label{5.32} 
A_\nu  =\frac{(-1)^{|\nu|}(|\sigma|+d+2|\nu|+1)_{N-|\nu|}}
   {(-N)_{|\nu|}\, N!}
 \prod_{j=1}^d \frac{(\sigma_j+a_j+\nu_j+1)_{\nu_j}(\sigma_j+1)_{\nu_j}\nu_j!}
       {(a_j+1)_{\nu_j}}.
\end{equation}
\end{Theorem}
 
These formulas are essentially contained in \cite{KM}. They can be deduced 
from \eqref{4.7} as in the proof of Theorem \ref{th5.2}.  Explicit 
biorthogonal (not mutually orthogonal) Hahn polynomials were also found in 
\cite{Tr1}. 
 

\subsubsection{} There is a nonempty set 
$S\subset\{1,2,\dots,d\}$ and $l_i+1 \in\{1,2,\dots,N\}$ for $i \in S$.
In this case the set $V$ is
\begin{equation} \label{5.33}
V^d_{N,S}=V^d_N\cap\{x:x_i\leq l_i\text{ for }i\in S\}.
\end{equation}
If $S=\{1,2,\dots,d\}$ we can also assume that $l_1+\cdots +l_d>N$, otherwise 
it becomes the parallelepiped case discussed in \thref{th5.3}. 

The weight $W(x)$ is again given by \eqref{5.29}, but $V^d_N$ is replaced by 
$V^d_{N,S}$. The corresponding polynomials are the same as in \thref{th5.4}, 
with the restriction $\nu\in V^d_{N,S}$.

\begin{Theorem}\label{th5.5}
For $\nu\in V^d_{N,S}$ the polynomials $Q_{\nu}(x;\sigma,N)$ defined by 
\eqref{5.30} satisfy the difference equation
$$ D Q_\nu = \lambda_{|\nu|}Q_\nu $$
and the orthogonal relation
\begin{equation*}
\sum_{x\in V} Q_\nu(x; \sigma,N) Q_\mu(x; \sigma, N)
\prod_{i=1}^{d}\binom{x_i+\sigma_i}{x_i} \binom{N-|x|+\sigma_{d+1}}{N-|x|}
 = A_\nu \delta_{\nu,\mu}, 
\end{equation*}
where $A_\nu$ is given by \eqref{5.32}.
\end{Theorem}

\begin{Remark}\label{re5.6}
Notice that \eqref{5.25} in the parallelepiped case $V^d_l$ gives a polynomial 
of total degree $|\nu|$. We can use the generators $G^d_l$ in \exref{ex3.1} to 
write the same polynomial using only the monomials $x^{\mu}$ with 
$\mu\in V^d_l$. Similarly, in the case of \thref{th5.5} we can use the
generators $\{(-x_i)_{l_i+1}:i\in S\}$ to express $Q_{\nu}$ in terms of the 
monomials $x^{\mu}$ with $\mu\in V^d_{N,S}$.
\end{Remark}

\begin{Remark}\label{re5.7}
In Theorems \ref{th5.3}-\ref{th5.5} we have 
$r_k=|\Lambda_{k}(V)|< \dim(\Pi^d_{k}/\Pi^d_{k-1})$ for some $k$'s 
and therefore the equation $Du(x)=\lambda u(x)$ will hold 
a priori only if we consider $u(x)$ as an element of $\R[V]$. However, one 
can show that the corresponding 
polynomials can be obtained as a limit from the polynomials in \thref{th5.2}.
This fact can be used to show that the equation $Du(x)=\lambda u(x)$ actually 
holds in the space $\R[x_1,x_2,\dots,x_d]$.
\end{Remark}


\subsection{Summary}\label{sse5.3}  
If $a= -1$ we have essentially one difference equation with coefficients
given in \eqref{5.2}. By specifying the free parameters, however, we end up 
with four different types of solutions given in Theorems \ref{th5.2} - 
\ref{th5.5}, respectively. 

As an example, let us consider the case $d=2$.  The weight functions and the
corresponding  sets on which they live are listed below:
\begin{itemize}
\item[(i)] $W(x) = \dfrac{(\sigma_1+1)_{x_1}(\sigma_2+1)_{x_2}}{x_1!x_2!}
       \dfrac{(\beta+1)_{|x|}}{(-|\sigma|-d-r+1)_{|x|}}$, \quad 
  $V = \N_0^2$;
\item[(ii)] $W(x)=\dfrac{(-l_1)_{x_1} (-l_2)_{x_2}}{x_1!x_2!}
       \dfrac{(\beta+1)_{|x|}}{(-l-r)_{|x|}}$, \quad $V= V^2_l$;
\item[(iii)] $W(x) = \dfrac{(\sigma_1+1)_{x_1}(\sigma_2+1)_{x_2}}{x_1!x_2!}
       \dfrac{(\sigma_3+1)_{N-|x|}}{(N-|x|)!}$, \quad  $V=V^2_N$;
\item[(iv)]  $W(x) =\dfrac{(-l_1)_{x_1}}{x_1!}\dfrac{(-l_2)_{x_2}}{x_2!} 
\dfrac{(\sigma_3+1)_{N-|x|}}{(N-|x|)!}$, 
                         \quad $V=V^2_{N,S}$, where $S=\{1,2\}$;
\item[(v)]  $W(x) =\dfrac{(-l_1)_{x_1}}{x_1!}\dfrac{(\sigma_2+1)_{x_2}}{x_2!} 
     \dfrac{(\sigma_3+1)_{N-|x|}}{(N-|x|)!}$, 
 \quad $V=V^2_{N,S}$, where $S=\{1\}$. 
\end{itemize}
In the last case, one can also exchange $x_1$ and $x_2$ to get another case. 
It should be mentioned that the equation \eqref{5.2} was considered in \
\cite{Xu2}, but only the case (iii) was identified there. 


\section{Multivariable case with linear eigenvalue: $a=0$ and $b=-1$} 
                                                               \label{se6}
When $a =0$ and $b = -1$, equations \eqref{3.15}-\eqref{3.18} become 
\begin{subequations}\label{6.1}
\begin{align}
&\alpha_{i,j}=-l_{i,j}x_j \quad \text{ for }i\neq j\label{6.1a}\\
&B_i=-x_i+s_i,\label{6.1b}\\
&\beta_i=\sum_{k\neq i}l_{i,k}x_k+(l_{i,i}-1)x_i+s_i\label{6.1c}\\
&\gamma_i=x_i\sum_{k=1}^dl_{k,i}.\label{6.1d}
\end{align}
\end{subequations}

First, notice that we cannot have $\gamma_i=0$. Indeed, if we assume that 
$\gamma_i=0$ then in order to have a self-adjoint operator we need 
$\beta_i=0$, which implies that $l_{i,k}=0$ for $k\neq i$. But then 
$\alpha_{i,k}=0$ for $k\neq i$ and therefore $\alpha_{k,i}=0$ for $k\neq i$, 
which simply means that the operator $D$ is independent of $x_i$ and 
$E_i^{\pm 1}$. Thus we have an operator acting 
in a $(d-1)$ dimensional space, trivially extended, by adding $x_i$ as 
explained in \reref{re3.6}.

Below we assume that $\gamma_i\neq 0$ for all $i$. Notice that $\gamma_i$ 
depends only on $x_i$, so in \eqref{2.12b} we can cancel $\gamma_i\gamma_j$ 
and we get 
\begin{equation}\label{6.2}
\beta_i(x-e_i)\beta_j(x-e_i-e_j)=
\beta_j(x-e_j)\beta_i(x-e_i-e_j).
\end{equation}
But $\beta_i(x-e_i)=\beta_i(x)+1-l_{i,i}$ and 
$\beta_i(x-e_i-e_j)=\beta_i(x)+1-l_{i,i}-l_{i,j}$ and plugging these in 
\eqref{6.2} we get
\begin{equation}\label{6.3}
l_{j,i}(\beta_i(x)+1-l_{i,i})=l_{i,j}(\beta_j(x)+1-l_{j,j}),
\end{equation}
for $i\neq j$. Comparing the coefficients of $x_k$ in the last formula we 
see that
\begin{align}
&l_{j,i}l_{i,k}=l_{i,j}l_{j,k} \text{ for }k\neq i,j\label{6.4}\\
&l_{j,i}l_{i,j}=l_{i,j}(l_{j,j}-1). \label{6.5}
\end{align}
These two equations have essentially two different types of solutions, which 
we discuss in two subsections.


\subsection{Case 1: $l_{i,j}\ne0$ for all $1\leq i\neq j\leq d$}\label{sse6.1}
 
\begin{Proposition} \label{pr6.1}
Assume that $l_{i,j}\neq 0$ for all $1\leq i\neq j\leq d$.  Then the most 
general solution of the compatibility conditions \eqref{2.12} is given by
\begin{subequations}\label{6.6}
\begin{align}
&\alpha_{i,j}=-l_{i}x_j \text{ for }i\neq j\label{6.6a}\\
&B_i=-x_i+l_is,\label{6.6b}\\
&\beta_i=l_i\left(\sum_{k=1}^dx_k +s\right)\label{6.6c}\\
&\gamma_i=x_i\left(\sum_{k=1}^dl_{k}+1\right),\label{6.6d}
\end{align}
\end{subequations}
where $s$ and $\{l_i\}_{i=1}^d$ are free parameters.
\end{Proposition}

\begin{proof} 
From equation \eqref{6.5} it follows that $l_{j,i}=l_{j,j}-1$. 
Denote $l_j:=l_{j,j}-1$. Then we have $l_{j,i}=l_j$ for all $i\neq j$. 
Using now \eqref{6.3} we see that $l_js_i=l_is_j$. Thus $s_i=l_is$,
which leads to formulas \eqref{6.6}.

Conversely, it is straightforward to see that if we define $\alpha_{i,j}$, 
$\beta_i$ and $\gamma_i$ as in \eqref{6.6a}, \eqref{6.6c} and \eqref{6.6d}, 
then the compatibility conditions in \coref{co2.3} are satisfied, i.e. the 
above formulas give the most general solution in the case $l_{i,j}\neq 0$.
\end{proof}

Below we determine the weight functions and the corresponding orthogonal 
polynomials. For every $i$ we have 
$$
\frac{W(x)}{W(x-e_i)}=\frac{l_i}{|l|+1}\frac{|x|+s-1}{x_i}
$$
and therefore 
$$
\frac{W(e_i)}{W(0)}=\frac{l_is}{|l|+1}>0,
$$
which shows that $l_i$ must have the same signs. There are two possible 
cases.


\subsubsection{ } $l_i=-p_i<0$ for all $i$. In this case, 
$$
\frac{W(x)}{W(x-e_i)}=\frac{p_i}{1-|p|}\frac{1-s-|x|}{x_i}.
$$
If the denominator does not vanish, then for $|x|$ large the second ratio 
will be negative and therefore we must have $1-|p|<0$, i.e. $|p|>1$. But 
then if we denote $c_i=p_i/(|p|-1)$ we will have $|c|>1$ and up to a constant 
factor the weight is $W(x)=(s)_{|x|}\prod_{i=1}^dc_i^{x_i}/x_i!$, which 
leads to divergent series.

Hence the only possibility here is $(-s)=N\in\N_0$. 
This forces $|p|<1$ and $V = V^d_N$. The difference 
operator is then 
\begin{equation} \label{6.7}
D=  \sum_{1\leq i,j\leq d}(\delta_{i,j}-p_i)x_j \fs_i\bs_j 
 + \sum_{i=1}^d (p_iN-x_i) \fs_i,
\end{equation}
and the eigenvalues are $\lambda_n = -n$. 
The orthogonal polynomials are the Krawtchouk polynomials on $V^d_N$. Recall
that Krawtchouk polynomial in one variable is denoted by $K_n(x;p,N)$.

\begin{Theorem}\label{th6.2}
Let $0< p_i < 1$, $1 \le i \le d$, and $|p|<1$. For $\nu\in\N_0^d$, 
$|\nu|\le N$, the polynomials
\begin{align}\label{6.8}
  K_\nu(x; p, N) = & \frac{(-1)^{|\nu|} }{(-N)_{|\nu|}}
    \prod_{j=1}^d  \frac{p_j^{\nu_j}} {(1-p_1-\cdots-p_j)^{\nu_j}}
       (-N+|X_{j-1}|+|\nu^{j+1}|)_{\nu_j} \\
  &\qquad\qquad  \times K_{\nu_j}\left(x_j; \tfrac{p_j}{1-p_1-\cdots-p_{j-1}}, 
         N-|X_{j-1}|- |\nu^{j+1}|\right) \notag
\end{align}
satisfy the difference equation
$$
     D \psi_\nu = \lambda_{|\nu|}\psi_\nu
$$
and the orthogonal relation
\begin{equation}\label{6.9}
 \sum_{|x| \le N} K_\nu(x;p,N)K_\mu(x;p,N) 
  \prod_{i=1}^{d+1} \frac{p_i^{x_i}}{x_i!}=\frac{(-1)^{|\nu|}}{(-N)_{|\nu|}N!}
     \prod_{j=1}^d \frac{ \nu_j! p_j^{\nu_j} } 
{(1-p_1-\cdots-p_j)^{\nu_j-\nu_{j+1}} } \delta_{\nu,\mu}, 
\end{equation}
where $x_{d+1}=N-|x|$, $p_{d+1}=1-|p|$ and $\nu_{d+1}=0$.
\end{Theorem}

In fact, these orthogonal polynomials can be considered as a limit of the 
Hahn polynomials \eqref{5.30} (see \cite{Tr2,Tr3}). Indeed, using the well
known relation
$$
   \lim_{t \to \infty} Q_n(x; p t, (1-p)t, N) = K_n(x; p, N)
$$
in one variable, it is not hard to see that 
$$
  \lim_{t \to \infty} Q_\nu(x; p_1 t, \ldots, p_d t, (1-p_1-\cdots - p_d)t, N)
    = K_\nu(x; p, N). 
$$
The orthogonality of $K_\nu(x; p, N)$ follows from \eqref{5.31} 
under the limit. We refer the reader to \cite{M} for other properties and 
applications of Krawtchouk polynomials.


\subsubsection{}  $l_i>0$ for all $i$.  In this case $s>0$. Denote 
$c_i=l_i/(1+|l|)$. Then  $c_i <1$ and $|c|< 1$,  the weight function is 
$$
W(x)=\frac{(s)_{|x|}}{x!}c^x =  (s)_{|x|}\prod_{i=1}^d \frac{c_i^{x_i}} {x_i!}.
$$
and $V = \N_0^d$. The difference operator takes the form 
\begin{equation} \label{6.10}
D=  \sum_{1\leq i,j\leq d}
\left(\delta_{i,j}+\frac{c_i}{1-|c|}\right)x_j \fs_i\bs_j 
+ \sum_{i=1}^d \left(- x_i + \frac{c_i}{1-|c|} s\right)\fs_i.
\end{equation}
The orthogonal polynomials are the Meixner polynomials on $\N_0^d$ but  they
are different from product Meixner polynomials. 
Recall that Meixner polynomial in one variable is denoted by $M_n(x;\beta,c)$.
We write $|C_j| = c_j  + c_{j+1}+ \cdots + c_d$ and define $C_{d+1} = 0$. 

\begin{Theorem}\label{th6.3}
Let $0< c_i < 1$, $1 \le i \le d$, and $|c|<1$. For $\nu \in \N_0^d$, 
the polynomials
\begin{align}\label{6.11}
  M_\nu(x; s, c) =   \prod_{j=1}^d  (\delta_j)_{\nu_j} 
               M_{\nu_j}\left(x_j; \delta_j , \frac{c_j}{1-|C_{j+1}|} \right)
\end{align}
where $\delta_j = s+|\nu^{j+1}|+|X_{j-1}|$, satisfy the difference equation
$$
     D \psi_\nu = \lambda_{|\nu|}\psi_\nu
$$
and the orthogonal relation
\begin{equation}\label{6.12}
  \sum_{x \in \N_0^d} M_\nu(x; s,c)M_\mu(x;s,c) 
    (s)_{|x|} \prod_{i=1}^{d} \frac{c_i^{x_i}}{x_i!} = 
\frac{(s)_{|\nu|}}{(1-|c|)^s}
\prod_{j=1}^d \nu_j! \left(\frac{c_j}{1-|C_{j+1}|}\right)^{-\nu_j}
          \delta_{\nu,\mu}. 
\end{equation}
\end{Theorem}

\begin{proof}
These relations can be derived inductively using the orthogonality 
\eqref{4.26} of the Meixner polynomials of one variable. It is clear from 
\eqref{4.24} that $M_{\nu}$ defined by \eqref{6.11} are indeed polynomials 
of $x$.
We can write the weight function as a product 
$$  
 W(x)=\frac{c_d^{x_d}(s+|X_{d-1}|)_{x_d}}{x_d!}W'(X_{d-1})
$$
where 
$$
W'( X_{d-1} ) = (s)_{|X_{d-1}|}\prod_{i=1}^{d-1} \frac{c_i^{x_i}}{x_i!}. 
$$
Using \eqref{4.26} the sum over $x_d$ in \eqref{6.12} becomes 
\begin{align*}
I_{\nu_d,\mu_d}&= (s+|X_{d-1}|)_{\nu_d}
        (s+|X_{d-1}|)_{\mu_d}
    \sum_{x_d=0}^{\infty}  M_{\nu_d}(x_d) M_{\mu_d}(x_d)
        \frac{c_d^{x_d}(s+|X_{d-1}|)_{x_d}}{x_d!}\\
&=\delta_{\nu_d,\mu_d}
\frac{c_d^{-\nu_d}\nu_d!}{(1-c_d)^s}
\frac{(s+|X_{d-1}|)_{\nu_d}}{(1-c_d)^{|X_{d-1}|}},
\end{align*}
where $M_m(x_d) = M_m(x_d,s+|X_{d-1}|, c_d)$. 

Combining this with $W'$ we obtain the new weight function
$$
W_{d-1}( X_{d-1} ) = 
\frac{c_d^{-\nu_d}\nu_d! (s)_{\nu_d}}{(1-c_d)^s}(s+\nu_d)_{|X_{d-1}|}
\prod_{i=1}^{d-1} \frac{1}{x_i!}\left(\frac{c_i}{1-c_d}\right)^{x_i}. 
$$

Apart from a constant multiple, the weight function $W_{d-1}$ has exactly the 
same structure of $W(x)$ with $s$ replaced by $s+\nu_d$, and $c_i$ replaced 
by $c_i/(1-c_d)$ for $i=1,2,\dots,d-1$. The proof now follows by induction. 
The constant $\langle M_{\nu},M_{\nu}\rangle$ is given by
$$\prod_{j=1}^d
\frac{\left(\frac{c_j}{1-|C_{j+1}|}\right)^{-\nu_j}
\nu_j! (s+|\nu^{j+1}|)_{\nu_j}}{{\left(1-\frac{c_j}{1-|C_{j+1}|}\right)^s}},
$$
which leads to \eqref{6.12}.
\end{proof}

\begin{Remark}\label{re6.4}
It is well known that the Meixner polynomials are the limit of the Hahn 
polynomials, 
$$
    \lim_{N \to \infty} Q_n(x;b-1,N\frac{1-c}{c},N)= M_n(x;b,c),
$$
There is an analogous relation between the Meixner polynomials \eqref{6.11} 
of several variables and the orthogonal polynomials on the parallelepiped
in \eqref{5.25}, at least when parameters $c$ have rational values. Indeed, 
let $\phi_\nu(x;\beta,r,l)$ denote the polynomials defined in \eqref{5.25}, 
which are orthogonal with respect to the weight function in \eqref{5.24}.  In 
these polynomials we set $\beta = s -1$ and 
$$
|L_{j-1}|+r=\left( \frac{1-|C_j|}{c_j} \right ) l_j, \qquad 1 \le j \le d.
$$
Using induction, it is easy to see that $l_j/c_j = l_{j-1}/c_{j-1}$, which 
then implies that 
$$
     l_j = (|l | +r)  c_j, \quad \hbox{or} \quad  l_j = \frac{c_j}{1-|c|} r,
     \qquad1 \le j \le d.
$$
If $c_j$ are rational numbers, we can choose $r = (1-|c|) N$, $N \in \N_0$, 
for certain $N$, so that $l_j$ are integers. Upon taking $r \to \infty$, which 
shows $l_j\to \infty$, it follows from \eqref{4.26} that 
\begin{equation} \label{6.14}
  \lim_{l \to \infty} \phi_\nu(x;s-1,r,l) = M_\nu(x;s,c).
\end{equation}
\end{Remark}


\subsection{Case 2: $l_{i,j}=0$ for some $1\leq i\neq j\leq d$} \label{sse6.2} 

Recall that $\alpha_{i,j}=0$ for $i\neq j$ if and only if $\alpha_{j,i}=0$. 
This essentially means that $l_{i,j}$ and $l_{j,i}$ are simultaneously zero 
or nonzero for $i\neq j$.

\begin{Lemma} \label{le6.5}
If $j\neq k$ and $l_{j,k}=0$ then for every $i\neq j,k$ we have 
$l_{i,j}=l_{j,i}=0$ or $l_{i,k}=l_{k,i}=0$.
\end{Lemma}
\begin{proof} Follows immediately from \eqref{6.4}.
\end{proof}

\begin{Theorem}\label{th6.6} 
Assume that $l_{i,j}=0$ for some $1\leq i\neq j\leq d$. 
Define $I=\{i\}\cup\{m:l_{m,i}\neq 0\}$ and $J=\{1,2,\dots,d\}\setminus I
\supset\{j\}$.
Then $D=D_I+D_J$, where $D_I$ is an admissible operator in the variables 
$\{x_m:m\in I\}$ and $D_J$ is an admissible operator in $\{x_k:k\in J\}$.
\end{Theorem}

\begin{proof} First we show that for $m\in I$ and $k\in J$ we have  
$l_{m,k}=l_{k,m}=0$. Indeed if $m=i$, then $l_{k,i}=0$ by the definition of 
$J$. If $m\neq i$ then $l_{m,i}\neq 0$, but $l_{k,i}=0$. The previous lemma  
shows that $l_{k,m}=l_{m,k}=0$.

Thus if $m\in I$ and $k\in J$ we see at once that:
\begin{itemize}
\item $\alpha_{m,k}=\alpha_{k,m}=0$
\item $\beta_m$ contains only the variables from $I$
\item $\beta_k$ contains only the variables from $J$.
\end{itemize}
The decomposition of $D$ follows immediately from the above observations.
\end{proof}

\begin{Remark}\label{re6.7}
\thref{th6.6} says that if $l_{i,j}=0$ 
for some $i\neq j$, then the operator $D$ splits as a sum of 2 operators of 
independent variables. Conversely let $D_I$ be an 
admissible operator in the variables $x'=\{x_m:m\in I\}$ and $D_J$ be an 
admissible operator in $x''=\{x_k:k\in J\}$ with $I\cap J=\emptyset$ and 
$I\cup J=\{1,2,\dots,d\}$. Let us denote by $p^{I}_{\mu}(x')$ and 
$p^{J}_{\nu}(x'')$ the polynomials satisfying
$$D_I(p^{I}_{\mu})=-|\mu|p^{I}_{\mu}\qquad \text{and}\qquad
D_J(p^{J}_{\nu})=-|\nu|p^{J}_{\nu}.$$
Then if put $D=D_I+D_J$ and $p_{\mu,\nu}(x)=p^{I}_{\mu}(x')p^{J}_{\nu}(x'')$,
where $x=(x',x'')\in\R^d$, we have 
$$Dp_{\mu,\nu}(x)=-(|\mu|+|\nu|)p_{\mu,\nu}(x).$$
\end{Remark}

Thus, in this case, the eigenfunctions of the difference operator are product 
of orthogonal polynomials of fewer variables, which satisfy difference 
equations of lower dimension with linear eigenvalues. For $d=1$ these are the  
polynomials of Charlier, Krawtchouk and Meixner given in Section \ref{sse4.2}. 
In higher dimensions, there are also the Krawtchouk polynomials of several 
variables given in \thref{th6.2} and the Meixner polynomials of several 
variables given in \thref{th6.3}.

Clearly there are many product polynomials of this type and the number 
increases drastically as the dimension grows.  As an example, we list all 
cases for $d =3$ below. To list the different types we use the abbreviation of 
C, K, M for Charlier, Meixner, and Krawtchouk polynomial of one variable, 
respectively, and use $\K2_2$ and $\M2_2$ to denote the Krawtchouk polynomials 
of two variables on $V^2_N$ and Meixner polynomials of two variables on 
$\N_0^2$. A product polynomial is denoted by its components. As an example, 
CCM stands for a product of the type Charlier-Charlier-Meixner.    
  
\begin{Example}[$d=3$]\label{ex6.8} 
There are sixteen product types, which we list according to 
their domains of orthogonality:
\begin{enumerate} 
\item $\N_0^3$: CCC, MMM, CCM, CMM, C$\M2_2$, M$\M2_2$.
\item $\N_0^2 \times V^1_N$: CCK, MMK, CMK, K$\M2_2$.
\item $\N_0 \times V^1_{N_1} \times V^1_{N_2}$,  CKK, MKK.
\item $V^1_{N_1} \times V^1_{N_2} \times V^1_{N_3}$:  KKK. 
\item $V^2_N \times \N_0$: $\K2_2$C, $\K2_2$M.
\item $V^2_{N_1} \times V^1_{N_2}$: $\K2_2$K. 
\end{enumerate} 
\end{Example}

\bigskip\noindent
{\bf Acknowledgments.} It is our pleasure to thank a referee for several 
thoughtful suggestions and for bringing the reference \cite{M} to our 
attention.

\end{document}